\documentclass[10.5pt,a4paper]{article}

\usepackage[leqno]{amsmath}
\usepackage{graphicx}
\usepackage{latexsym}
\usepackage{amsmath,amsfonts,amssymb,amsthm,euscript,makeidx,color,mathrsfs}
\usepackage{enumerate}

\usepackage[colorlinks,linkcolor=blue,anchorcolor=green,citecolor=red]{hyperref}

\def\sqr#1#2{{\vcenter{\vbox{\hrule height.#2pt
              \hbox{\vrule width.#2pt height#1pt \kern#1pt \vrule width.#2pt}
              \hrule height.#2pt}}}}
\def\signed #1{{\unskip\nobreak\hfil\penalty50
              \hskip2em\hbox{}\nobreak\hfil#1
              \parfillskip=0pt \finalhyphendemerits=0 \par}}
\def\endpf{\signed {$\sqr69$}}

\def\5n{\negthinspace \negthinspace \negthinspace \negthinspace \negthinspace }
\def\4n{\negthinspace \negthinspace \negthinspace \negthinspace }
\def\3n{\negthinspace \negthinspace \negthinspace }
\def\2n{\negthinspace \negthinspace }
\def\1n{\negthinspace }

\def\dbE{\mathbb{E}}
\def\dbF{\mathbb{F}}

\def\dbH{\mathbb{H}}

\def\dbP{\mathbb{P}}

\def\dbR{\mathbb{R}}
\def\dbS{\mathbb{S}}

\def\sC{\mathscr{C}}

\def\ds{\displaystyle}

\def\ns{\noalign{\ss}}
\def\a{\alpha}
\def\b{\beta}

\def\d{\delta}

\def\z{\zeta}

\def\si{\sigma}
\def\t{\tau}
\def\f{\varphi}
\def\th{\theta}

\def\i{\infty}
\def\G{\Gamma}
\def\D{\Delta}
\def\Th{\Theta}
\def\L{\Lambda}
\def\Si{\Sigma}
\def\F{\Phi}
\def\Om{\Omega}
\def\cF{{\cal F}}

\def\cR{{\cal R}}

\def\cU{{\cal U}}

\def\BJ{{\bf J}}

\def\BX{{\bf X}}
\def\BY{{\bf Y}}

\def\no{\noindent}

\def\ss{\smallskip}
\def\ms{\medskip}

\def\q{\quad}
\def\qq{\qquad}
\def\hb{\hbox}

\def\Ra{\mathop{\Rightarrow}}

\def\lt{\left}
\def\rt{\right}
\def\lan{\langle}
\def\ran{\rangle}
\def\llan{\left\langle}
\def\rran{\right\rangle}
\def\blan{\big\langle}
\def\bran{\big\rangle}

\def\rf{\eqref}

\def\h{\widehat}
\def\wt{\widetilde}

\def\cd{\cdot}

\def\ae{\hbox{\rm a.e.}}
\def\as{\hbox{\rm a.s.}}

\def\tr{\hbox{\rm tr$\,$}}

\def\deq{\triangleq}
\def\les{\leqslant}
\def\ges{\geqslant}

\def\({\Big (}
\def\){\Big )}
\def\[{\Big[}
\def\]{\Big]}
\def\pr{{\partial}}
\def\bde{\begin{definition}\label}
\def\ede{\end{definition}}
\def\be{\begin{equation}}
\def\bel{\begin{equation}\label}
\def\ee{\end{equation}}
\def\bt{\begin{theorem}\label}
\def\et{\end{theorem}}
\def\bc{\begin{corollary}\label}
\def\ec{\end{corollary}}
\def\bl{\begin{lemma}\label}
\def\el{\end{lemma}}
\def\bp{\begin{proposition}\label}
\def\ep{\end{proposition}}
\def\bas{\begin{assumption}\label}
\def\eas{\end{assumption}}
\def\br{\begin{remark}\label}
\def\er{\end{remark}}
\def\bex{\begin{example}\label}
\def\ex{\end{example}}
\def\ba{\begin{array}}
\def\ea{\end{array}}
\def\ed{\end{document}}

\def\square#1{\vbox{\hrule\hbox{\vrule height#1%
     \kern#1\vrule}\hrule}}
\def\rectangle#1#2{\vbox{\hrule\hbox{\vrule height#1%
     \kern#2\vrule}\hrule}}

\font\tenbb=msbm10 \font\sevenbb=msbm7 \font\fivebb=msbm5

\newfam\bbfam
\scriptscriptfont\bbfam=\fivebb \textfont\bbfam=\tenbb
\scriptfont\bbfam=\sevenbb

\newtheorem{theorem}{\hskip 1.3em Theorem}[section]
\newtheorem{definition}[theorem]{\hskip 1.3em Definition}
\newtheorem{proposition}[theorem]{\hskip 1.3em Proposition}
\newtheorem{corollary}[theorem]{\hskip 1.3em Corollary}
\newtheorem{lemma}[theorem]{\hskip 1.3em Lemma}
\newtheorem{remark}[theorem]{\hskip 1.3em Remark}
\newtheorem{example}[theorem]{\hskip 1.3em Example}

\newtheorem{assumption}[theorem]{\hskip 1.3em Assumption}

\makeatletter
   
   \@addtoreset{equation}{section}
\makeatother

\begin{document}

\title{\bf Mean-Field Stochastic Linear Quadratic Optimal Control Problems: Closed-Loop
Solvability}
\author{Xun Li\thanks{Department of Applied Mathematics, The Hong Kong Polytechnic University,
Hong Kong, China (malixun@polyu.edu.hk). This author was partially supported by Hong Kong RGC
under grants 519913, 15209614 and 15224215.}\ , \ Jingrui Sun\thanks{Department of Applied Mathematics, The Hong Kong Polytechnic University, Hong Kong, China
(sjr@mail.ustc.edu.cn). This author was partially supported by the National Natural Science
Foundation of China (11401556) and the Fundamental Research Funds for the Central Universities
(WK 2040000012).}\ , \ and \ Jiongmin Yong\thanks{Department of Mathematics, University of
Central Florida, Orlando, FL 32816, USA (jiongmin.yong@ucf.edu). This author was partially
supported by NSF DMS-1406776.}}
\maketitle

\no\bf Abstract: \rm An optimal control problem is studied for a linear mean-field stochastic
differential equation with a quadratic cost functional. The coefficients and the weighting
matrices in the cost functional are all assumed to be deterministic. Closed-loop strategies
are introduced, which require to be independent of initial states; and such a nature makes
it very useful and convenient in applications. In this paper, the existence of an optimal
closed-loop strategy for the system (also called the closed-loop solvability of the problem)
is characterized by the existence of a regular solution to the coupled two (generalized)
Riccati equations, together with some constraints on the adapted solution to a linear backward
stochastic differential equation and a linear terminal value problem of an ordinary differential
equation.

\ms

\no\bf Key words: \rm mean-field stochastic differential equation, linear quadratic optimal
control, Riccati equation, regular solution, closed-loop solvability

\ms

\no\bf AMS subject classifications. \rm 49N10, 49N35, 93E20

\section{Introduction}


\ms

Let $(\Om,\cF,\dbF,\dbP)$ be a complete filtered probability space on which a standard
one-dimensional Brownian motion $W=\{W(t); 0\les t<\i\}$ is defined,
where $\dbF=\{\cF_t\}_{t\ges0}$ is the natural filtration of $W$ augmented by all the
$\dbP$-null sets in $\cF$. Consider the following controlled linear mean-field stochastic differential equation (MF-SDE, for short) on a finite time horizon $[t,T]$:
\bel{state}\left\{\2n\ba{ll}
\ds dX(s)=\Big\{A(s)X(s)+\bar A(s)\dbE[X(s)]+B(s)u(s)+\bar B(s)\dbE[u(s)]+b(s)\Big\}ds\\
\ns\ds\qq\qq~
+\Big\{C(s)X(s)+\bar C(s)\dbE[X(s)]+D(s)u(s)+\bar D(s)\dbE[u(s)]+\si(s)\Big\}dW(s),\qq s\in[t,T],\\
\ns\ds X(t)=\xi,
\ea\right.\ee
where $A(\cd)$, $\bar A(\cd)$, $B(\cd)$, $\bar B(\cd)$, $C(\cd)$, $\bar C(\cd)$, $D(\cd)$,
$\bar D(\cd)$ are given deterministic matrix-valued functions; $b(\cd)$, $\si(\cd)$ are
vector-valued $\dbF$-progressively measurable processes and $\xi$ is an $\cF_t$-measurable
random vector. In the above, $u(\cd)$ is the {\it control process} and $X(\cd)$ is the
corresponding {\it state process} with {\it initial pair} $(t,\xi)$. For any $t\in[0,T)$,
we define
$$\cU[t,T]=\lt\{u:[t,T]\times\Om\to\dbR^m\bigm|u(\cd)\hb{ is $\dbF$-progressively measurable, }
\dbE\int_t^T|u(s)|^2ds<\i\rt\}.$$
Any $u(\cd)\in\cU[t,T]$ is called an {\it admissible control} (on $[t,T]$). Under some mild conditions,
for any initial pair $(t,\xi)$ with $\xi\in L^2_{\cF_t}(\Om;\dbR^n)$ (the set of all $\cF_t$-measurable,
square-integrable $\dbR^n$-valued processes), and any admissible control $u(\cd)\in\cU[t,T]$, \rf{state}
admits a unique square-integrable solution $X(\cd)\equiv X(\cd\,;t,\xi,u(\cd))$.
Now we introduce the following cost functional:
\bel{cost}\ba{ll}
\ds J(t,\xi;u(\cd))\deq\dbE\Bigg\{\lan GX(T),X(T)\ran+2\lan g,X(T)\ran
+\blan\bar G\dbE[X(T)],\dbE[X(T)]\bran+2\lan\bar g,\dbE[X(T)]\ran\\
\ns\ds\qq\qq\,+\int_t^T\lt[\llan\begin{pmatrix}Q(s)&S(s)^\top\\S(s)&R(s)\end{pmatrix}
                                \begin{pmatrix}X(s)\\ u(s)\end{pmatrix},
                                \begin{pmatrix}X(s)\\u(s)\end{pmatrix}\rran
+2\llan\begin{pmatrix}q(s)\\ \rho(s)\end{pmatrix},\begin{pmatrix}X(s)\\ u(s)\end{pmatrix}\rran\rt]ds\\
\ns\ds\qq\qq\,+\int_t^T\lt[\llan\begin{pmatrix}\bar Q(s)&\bar S(s)^\top\\\bar S(s)&\bar R(s)\end{pmatrix}
                                \begin{pmatrix}\dbE[X(s)]\\ \dbE[u(s)]\end{pmatrix},
                                \begin{pmatrix}\dbE[X(s)]\\ \dbE[u(s)]\end{pmatrix}\rran
+2\llan\begin{pmatrix}\bar q(s)\\ \bar\rho(s)\end{pmatrix},
       \begin{pmatrix}\dbE[X(s)]\\ \dbE[u(s)]\end{pmatrix}\rran\rt] ds\Bigg\},
\ea\ee
where $G$, $\bar G$ are symmetric matrices and $Q(\cd)$, $\bar Q(\cd)$, $S(\cd)$, $\bar S(\cd)$,
$R(\cd)$, $\bar R(\cd)$ are deterministic matrix-valued functions with $Q(\cd)$, $\bar Q(\cd)$,
$R(\cd)$, and $\bar R(\cd)$ being symmetric; $g$ is an $\cF_T$-measurable random vector and
$\bar g$ is a deterministic vector; $q(\cd)$, $\rho(\cd)$ are vector-valued $\dbF$-progressively
measurable processes and $\bar q(\cd)$, $\bar \rho(\cd)$ are vector-valued deterministic functions.
Our mean-field stochastic linear quadratic (LQ, for short) optimal control problem can be stated
as follows:

\ms

\bf Problem (MF-LQ). \rm For any given initial pair $(t,\xi)\in[0,T)\times L^2_{\cF_t}(\Om;\dbR^n)$,
find a $u^*(\cd)\in\cU[t,T]$ such that
\bel{optim}J(t,\xi;u^*(\cd))=\inf_{u(\cd)\in\cU[t,T]}J(t,\xi;u(\cd))\deq V(t,\xi).\ee

Any $u^*(\cd)\in\cU[t,T]$ satisfying \rf{optim} is called an {\it optimal open-loop control}
of Problem (MF-LQ) for the initial pair $(t,\xi)$, and the corresponding
$X^*(\cd)\equiv X(\cd\,; t,\xi,u^*(\cd))$ is called an {\it optimal open-loop state process}.
The function $V(\cd\,,\cd)$ is called the {\it value function} of Problem (MF-LQ).
In the special case where $b(\cd)$, $\si(\cd)$, $g$, $\bar g$, $q(\cd)$, $\bar q(\cd)$, $\rho(\cd)$,
and  $\bar\rho(\cd)$ vanish, we denote the corresponding mean-field LQ problem, cost functional,
and value function by Problem (MF-LQ)$^0$, $J^0(t,\xi;u(\cd))$, and $V^0(t,\xi)$, respectively.

\ms

The theory of MF-SDEs can be traced back to Kac who presented a stochastic toy model for the Vlasov
kinetic equation of plasma in \cite{Kac 1956} which leads to the so-called McKean-Vlasov
stochastic differential equation. Since then, researches on the related topics and their
applications have become a notable and serious endeavor among researchers in applied
probability and optimal stochastic controls, including financial engineering.
See, for examples, McKean \cite{McKean 1966},
Buckdahn--Djehiche--Li--Peng \cite{Buckdahn-Djehiche-Li-Peng 2009},
Buckdahn--Li--Peng \cite{Buckdahn-Li-Peng 2009}, Andersson--Djenhiche \cite{Andersson-Djehiche 2011},
Buckdahn--Djehiche--Li \cite{Buckdahn-Djehiche-Li 2011},
Meyer-Brandis--{\O}ksendal--Zhou \cite{Meyer-Oksendal-Zhou 2011},
Yong \cite{Yong 2013}, Elliott--Li--Ni \cite{Elliott-Li-Ni 2013},
Cui--Li--Li \cite{Cui-Li-Li 2014}, Huang--Li--Wang \cite{Huang-Li-Wang 2015}, Huang--Li--Yong \cite{Huang-Li-Yong 2015}. Note that when the mean-field part is absent, Problem (MF-LQ) is reduced
to the classical stochastic LQ optimal control problem. For relevant results and historic remarks on this subject, the reader is further referred to, for examples,
\cite{Wonham 1968,Chen-Li-Zhou 1998,Chen-Yong 2000,Ait Rami-Moore-Zhou 2001,Tang 2003}
and the book of Yong--Zhou \cite{Yong-Zhou 1999}.

\ms

Recently, Sun and Yong introduced the notions of open-loop and closed-loop solvabilities for stochastic LQ problems \cite{Sun-Yong 2014}. See also \cite{Sun-Yong-Zhang 2016, Sun-Li-Yong}. It turns out that these two notions are essentially different for stochastic LQ problem on finite time horizon. Roughly speaking, the open-loop solvability is equivalent to the solvability of a forward-backward stochastic differential equation (FBSDE, for short), and the closed-loop solvability is equivalent to the existence of a regular solution to a Riccati equation. Open-loop solvability was studied for Problem (MF-LQ) in \cite{Sun}. The current work is therefore a continuation of the above-mentioned works.

\ms

The rest of the paper is organized as follows. In Section 2 we give some preliminaries,
carefully explain the closed-loop strategies, and introduce the regular solution to the
generalized Riccati equations. Section 3 is devoted to the necessary conditions for the
existence of an optimal closed-loop strategy. In Section 4, we present our main result,
in which the closed-loop solvability of the mean-field LQ problems is characterized.
Finally, some concluding remarks are given in Section 5.

\section{Preliminaries}

We begin with some notation that will be used throughout the paper:
$$\ba{lll}
\hb{$\dbR^{n\times m}$: the Euclidean space of all $n\times m$ real matrices;
$\dbR^n=\dbR^{n\times 1}$ and $\dbR=\dbR^1$}.\\
\hb{$\dbS^n$: the space of all symmetric $n\times n$ real matrices}.\\
%
%
%
\hb{$M^\top$: the transpose of a matrix $M$}.\\
\hb{$M^\dag$: the Moore--Penrose pseudoinverse of a matrix $M$ \cite{Penrose 1955}}.\\
\hb{$\tr(M)$: the sum of diagonal elements of a square matrix $M$}.\\
\hb{$\lan\cd\,,\cd\ran$: the inner product on a Euclidean space given by $\lan M,N\ran\mapsto\tr(M^\top N)$.}\\
\hb{$|M|\deq\sqrt{\tr(MM^\top)}$: the Frobenius nrom of a matrix $M$}.\\
\hb{$\cR(M)$: the range of a matrix $M$}.\ea$$
For $M,N\in\dbS^n$, we use the notation $M\ges N$ (respectively, $M>N$) to indicate
that $M-N$ is positive semi-definite (respectively, positive definite).
For any $\dbS^n$-valued measurable function $F$ on $[t,T]$, we write
$$\ba{llll}
\ds   F\ges0 & \Longleftrightarrow &\ds F(s)\ges0,    \q& \ae~s\in[t,T],\\
\ns\ds   F>0 & \Longleftrightarrow &\ds F(s)>0,       \q& \ae~s\in[t,T],\\
\ns\ds F\gg0 & \Longleftrightarrow &\ds F(s)\ges\d I, \q& \ae~s\in[t,T],\hb{ for some }\d>0.
\ea$$
Let $f(\cd)$ be a function from $\dbR^{n\times m}$ into $\dbR$. Recall that the gradient of
$f$ at $X=(x_{ij})$, denoted by ${\pr f(X)\over\pr X}$, is an $n\times m$ matrix whose $(i,j)$-th
entry is given by $\big[{\pr f(X)\over\pr X}\big]_{ij}={\pr f(X)\over\pr x_{ij}}$.
For matrices $L$, $M$, and $N$ of proper dimensions, the following formulae hold:
$$\ba{lll}
\no\ms&\ds{\pr\over\pr X}\tr(LXM)=L^\top\1n M^\top,\q
&\ds{\pr\over\pr X}\tr(X^\top\1n LXM)=LXM+L^\top\1n XM^\top,\\
&\ds{\pr\over\pr X}\tr(LX^\top\1n M)=ML,
\q &\ds{\pr\over\pr X}\tr(LXMX^\top\1n N)=L^\top\1nN^\top\1nXM+NLXM. \ea$$

\ss

For a Euclidean space $\dbH$, let $L^p(t,T;\dbH)$ $(1\les p\les\i)$ be the space of all
$\dbH$-valued functions that are $L^p$-integrable on $[t,T]$, and let $C([t,T];\dbH)$ be
the space of all $\dbH$-valued continuous functions on $[t,T]$. We denote
$$\ba{ll}
\ds L^2_{\cF_t}(\Om;\dbH)=\Big\{\xi:\Om\to\dbH\bigm|\xi\hb{ is $\cF_t$-measurable, }
\dbE|\xi|^2<\i\Big\},\\
\ns\ds L_\dbF^2(t,T;\dbH)=\Big\{\f:[t,T]\times\Om\to\dbH\bigm|\f(\cd)\hb{ is
$\dbF$-progressively measurable, }\dbE\int^T_t|\f(s)|^2ds<\i\Big\},\\
\ns\ds L_\dbF^2(\Om;C([t,T];\dbH))=\Big\{\f:[t,T]\times\Om\to\dbH\bigm|\f(\cd)\hb{ is
$\dbF$-adapted, continuous, }\dbE\(\sup_{s\in[t,T]}|\f(s)|^2\)<\i\Big\},\\
\ns\ds L^2_\dbF(\Om;L^1(t,T;\dbH))=\Big\{\f:[t,T]\times\Om\to\dbH\bigm|\f(\cd)\hb{ is
$\dbF$-progressively measurable, }\dbE\(\int_t^T|\f(s)|ds\)^2<\i\Big\}.
\ea$$
The following assumptions will be in force throughout this paper.

\ms

{\bf(H1)} The coefficients of the state equation satisfy the following:
$$\left\{\2n\ba{ll}
\ds A(\cd),\bar{A}(\cd)\in L^1(0,T;\dbR^{n\times n}),
\q  B(\cd),\bar{B}(\cd)\in L^2(0,T;\dbR^{n\times m}),
\q  b(\cd)\in L^2_\dbF(\Om;L^1(0,T;\dbR^n)),\\
\ns\ds C(\cd),\bar C(\cd)\in L^2(0,T;\dbR^{n\times n}),
\q     D(\cd),\bar D(\cd)\in L^{\i}(0,T;\dbR^{n\times m}),
\q     \si(\cd)\in L^2_\dbF(0,T;\dbR^n).
\ea\right.$$

{\bf(H2)} The weighting coefficients in the cost functional satisfy the following:
$$\left\{\2n\ba{ll}
\ds Q(\cd), \bar Q(\cd)\in L^1(0,T;\dbS^n),
\q  S(\cd), \bar S(\cd)\in L^2(0,T;\dbR^{m\times n}),
\q  R(\cd), \bar R(\cd)\in L^\i(0,T;\dbS^m),\\
\ns\ds g\in L^2_{\cF_T}(\Om;\dbR^n),
\q     q(\cd)\in L^2_\dbF(\Om;L^1(0,T;\dbR^n)),
\q     \rho(\cd)\in L_\dbF^2(0,T;\dbR^m),\\
\ns\ds\bar g\in \dbR^n,\q \bar q(\cd)\in L^1(0,T;\dbR^n),
\q    \bar\rho(\cd)\in L^2(0,T;\dbR^m),\q G, \bar G\in\dbS^n.
\ea\right.$$

\ss

A standard argument using the contraction mapping theorem shows that under (H1),
for any initial pair $(t,\xi)\in[0,T)\times L^2_{\cF_t}(\Om;\dbR^n)$ and any
admissible control $u(\cd)\in\cU[t,T]$, \rf{state} admits a unique (strong) solution
$X(\cd)\equiv X(\cd\,; t,\xi,u(\cd))\in L_\dbF^2(\Om;C([t,T];\dbR^n))$.
Hence, under (H1)--(H2), the cost functional \rf{cost} is well-defined,
and Problem (MF-LQ) makes sense.

\ms

Let us now recall the notion of open-loop solvability from Sun \cite{Sun},
which was inspired by \cite{Sun-Yong 2014}.

\bde{def-open}\rm Problem (MF-LQ) is said to be ({\it uniquely}) {\it open-loop solvable
at initial pair $(t,\xi)\in[0,T]\times L^2_{\cF_t}(\Om;\dbR^n)$} if there exists a (unique)
$u^*(\cd)\in\cU[t,T]$ satisfying \rf{optim}.
Problem (MF-LQ) is said to be ({\it uniquely}) {\it open-loop solvable at $t$} if for any
$\xi\in L^2_{\cF_t}(\Om;\dbR^n)$, there exists a (unique) $u^*(\cd)\in\cU[t,T]$ satisfying \rf{optim},
and Problem (MF-LQ) is said to be ({\it uniquely}) {\it open-loop solvable on $[t,T)$}
if it is (uniquely) open-loop solvable at all $s\in[t,T)$.
\ede

Next, inspired by \cite{Sun-Yong 2014}, we introduce the following definition.

\bde{def-closed}\rm (i) Let
$$\sC[t,T]=L^2(t,T;\dbR^{m\times n})\times L^2(t,T;\dbR^{m\times n})\times\cU[t,T].$$
Any triple $(\Th(\cd),\bar\Th(\cd),v(\cd))\in\sC[t,T]$ is called a {\it closed-loop strategy}
of Problem (MF-LQ) on $[t,T]$.

\ms

(ii) For any $(\Th(\cd),\bar\Th(\cd),v(\cd))\in\sC[t,T]$ and $\xi\in L^2_{\cF_t}(\Om;\dbR^n)$,
let $X(\cd)\equiv X(\cd\,;t,\xi,\Th(\cd),\bar\Th(\cd),v(\cd))$ be the solution to the following
{\it closed-loop system}:
\bel{Feb12-01}\left\{\2n\ba{ll}
\ds dX(s)=\Big\{\big[A(s)+B(s)\Th(s)\big]X(s)+\big\{\bar A(s)+B(s)\bar\Th(s)
+\bar B(s)\big[\Th(s)+\bar\Th(s)\big]\big\}\dbE[X(s)]\\
\ns\ds\qq\qq\qq~+B(s)v(s)+\bar B(s)\dbE[v(s)]+b(s)\Big\}ds\\
\ns\ds\qq\q~~~+\Big\{\big[C(s)+D(s)\Th(s)\big]X(s)+\big\{\bar C(s)+D(s)\bar\Th(s)
+\bar D(s)\big[\Th(s)+\bar\Th(s)\big]\big\}\dbE[X(s)]\\
\ns\ds\qq\qq\qq~+D(s)v(s)+\bar D(s)\dbE[v(s)]+\si(s)\Big\}dW(s),\qq s\in[t,T],\\
\ns\ds X(t)=\xi,
\ea\right.\ee
and let
$$u(s)=\Th(s)X(s)+\bar\Th(s)\dbE[X(s)]+v(s),\qq s\in[t,T].$$
Then $(X(\cd),u(\cd))$ is called the {\it outcome pair} of $(\Th(\cd),\bar\Th(\cd),v(\cd))$
on $[t,T]$ corresponding to the initial state $\xi$; $X(\cd)$ and $u(\cd)$ are called the
corresponding {\it closed-loop state process} and {\it closed-loop outcome control process}, respectively.
\ede

Note that if $(X(\cd),u(\cd))$ is the outcome pair of a closed-loop strategy
$(\Th(\cd),\bar\Th(\cd),v(\cd))\in\sC[t,T]$ corresponding to some $\xi\in L^2_{\cF_t}(\Om;\dbR^n)$,
then $(X(\cd),u(\cd))\in L^2_\dbF(\Om;C([t,T];\dbR^n)\times\cU[t,T]$ is actually a state-control pair
of the state equation \rf{state}. In fact, \rf{Feb12-01} is equivalent to the following:
$$\left\{\2n\ba{ll}
\ds dX(s)=\Big\{A(s)X(s)+\bar A(s)\dbE[X(s)]+B(s)\big\{\Th(s)X(s)+\bar\Th(s)\dbE[X(s)]+v(s)\big\}\\
\ns\ds\qq\qq\qq~
+\bar B(s)\dbE\big\{\Th(s)X(s)+\bar\Th(s)\dbE[X(s)]+v(s)\big\}+b(s)\Big\}ds\\
\ns\ds\qq\q~~~+\Big\{C(s)X(s)+\bar C(s)\dbE[X(s)]+D(s)\big\{\Th(s)X(s)+\bar\Th(s)\dbE[X(s)]+v(s)\big\}\\
\ns\ds\qq\qq\qq~
+\bar D(s)\dbE\big\{\Th(s)X(s)+\bar\Th(s)\dbE[X(s)]+v(s)\big\}+\si(s)\Big\}dW(s),\qq s\in[t,T],\\
\ns\ds X(t)=\xi.
\ea\right.$$
Therefore, for any $(\Th(\cd),\bar\Th(\cd),v(\cd))\in\sC[t,T]$, the meaning of
$J(t,\xi;\Th(\cd)X(\cd)+\bar\Th(\cd)+v(\cd))$ is clear. We point out that a closed-loop strategy $(\Th(\cd),\bar\Th(\cd),v(\cd))\in\sC[t,T]$ is not related to any initial state
$\xi\in L^2_{\cF_t}(\Om;\dbR^n)$, whereas, an outcome pair $(X(\cd),u(\cd))$ of a closed-loop strategy $(\Th(\cd),\bar\Th(\cd),v(\cd))$ depends not only on the closed-loop strategy, but also on the initial
state $\xi$. Hence, we should carefully distinguish the closed-loop strategy and the corresponding
outcome control. Now, we are ready to introduce the following notion.

\bde{optimal closed-loop strategy} \rm
(i) A closed-loop strategy  $(\Th^*(\cd),\bar\Th^*(\cd),v^*(\cd))\in\sC[t,T]$ is said to be
{\it optimal} on $[t,T]$ if
\bel{strong}\ba{ll}
\ds J(t,\xi;\Th^*(\cd)X^*(\cd)+\bar\Th^*(\cd)\dbE[X^*(\cd)]+v^*(\cd))
\les J(t,\xi;\Th(\cd)X(\cd)+\bar\Th(\cd)\dbE[X(\cd)]+v(\cd)),\\
\ns\ds\qq\qq\qq\qq\qq\qq\qq\qq\1n\,\forall(\Th(\cd),\bar\Th(\cd),v(\cd))\in\sC[t,T],
~\forall\xi\in L^2_{\cF_t}(\Om;\dbR^n),
\ea\ee
where $X^*(\cd)$ and $X(\cd)$ are the closed-loop state processes corresponding to
$(\Th^*(\cd),\bar\Th^*(\cd),v^*(\cd),\xi)$ and $(\Th(\cd)$, $\bar\Th(\cd),v(\cd),\xi)$,
respectively. If an optimal closed-loop strategy (uniquely) exists on $[t,T]$,
Problem (MF-LQ) is said to be ({\it uniquely}) {\it closed-loop solvable} on $[t,T]$.

\ms

(ii) A closed-loop strategy $(\Th^*(\cd),\bar\Th^*(\cd),v^*(\cd))\in\sC[t,T]$ is said to be
{\it weakly optimal} on $[t,T]$ if \rf{strong} holds only for $\xi=x\in\dbR^n$. If a weakly
optimal closed-loop strategy (uniquely) exists on $[t,T]$, Problem (MF-LQ) is said to be
({\it uniquely}) {\it weakly closed-loop solvable} on $[t,T]$.
\ede

Similar to \cite{Sun-Yong 2014}, we have the following proposition.

\bp{equivalence}\sl Let {\rm(H1)--(H2)} hold and let $(\Th^*(\cd),\bar\Th^*(\cd),v^*(\cd))\in\sC[t,T]$.
Then the following statements are equivalent:
\begin{enumerate}[~~\,\rm(i)]
\item $(\Th^*(\cd),\bar\Th^*(\cd),v^*(\cd))$ is an optimal closed-loop strategy of
Problem {\rm(MF-LQ)} on $[t,T]$;
\item The following holds:
$$\ba{ll}
\ds J(t,\xi;\Th^*(\cd)X^*(\cd)+\bar\Th^*(\cd)\dbE[X^*(\cd)]+v^*(\cd))
\les J(t,\xi;\Th^*(\cd)X(\cd)+\bar\Th^*(\cd)\dbE[X(\cd)]+v(\cd)),\\
\ns\ds\qq\qq\qq\qq\qq\qq\qq\qq\qq\qq\qq~~\forall(\xi,v(\cd))\in L^2_{\cF_t}(\Om;\dbR^n)\times\cU[t,T],
\ea$$
where $X^*(\cd)$ and $X(\cd)$ are the closed-loop state processes corresponding to
$(\Th^*(\cd),\bar\Th^*(\cd),v^*(\cd),\xi)$ and $(\Th^*(\cd),\bar\Th^*(\cd),v(\cd),\xi)$,
respectively;
\item The following holds:
\bel{strong**} J(t,\xi;\Th^*(\cd)X^*(\cd)+\bar\Th^*(\cd)\dbE[X^*(\cd)]+v^*(\cd))\les
J(t,\xi;u(\cd)),\q\forall(\xi,u(\cd))\in L^2_{\cF_t}(\Om;\dbR^n)\times\cU[t,T],\ee
where $X^*(\cd)$ is the closed-loop state process corresponding to
$(\Th^*(\cd),\bar\Th^*(\cd),v^*(\cd))$ and $\xi\in L^2_{\cF_t}(\Om;\dbR^n)$.
\end{enumerate}
\ep

\it Proof. \rm (i) $\Ra$ (ii) is trivial, by taking $\Th(\cd)=\Th^*(\cd)$
and $\bar\Th(\cd)=\bar\Th^*(\cd)$ in \rf{strong}.

\ms

(ii) $\Ra$ (iii). For any $\xi\in L^2_{\cF_t}(\Om;\dbR^n)$ and $u(\cd)\in\cU[t,T]$,
let $X(\cd)$ be the solution to the following:
\bel{12.25-1}\left\{\2n\ba{ll}
\ds dX(s)=\Big\{A(s)X(s)+\bar A(s)\dbE[X(s)]+B(s)u(s)+\bar B(s)\dbE[u(s)]+b(s)\Big\}ds\\
\ns\ds\qq\qq~+\Big\{C(s)X(s)+\bar C(s)\dbE[X(s)]+D(s)u(s)+\bar D(s)\dbE[u(s)]+\si(s)\Big\}dW(s),\qq s\in[t,T],\\
\ns\ds X(t)=\xi,\ea\right.\ee
and set
$$v(\cd)\deq u(\cd)-\Th^*(\cd)X(\cd)-\bar\Th^*(\cd)\dbE[X(\cd)]\in\cU[t,T].$$
Then $X(\cd)$ is also the solution to the following MF-SDE:
$$\left\{\2n\ba{ll}
\ds dX(s)=\Big\{A(s)X(s)+\bar A(s)\dbE[X(s)]+B(s)\big\{\Th^*(s)X(s)+\bar\Th^*(s)\dbE[X(s)]+v(s)\big\}\\
\ns\ds\qq\qq\qq~+\bar B(s)\dbE\big\{\Th^*(s)X(s)+\bar\Th^*(s)\dbE[X(s)]+v(s)\big\}+b(s)\Big\}ds\\
\ns\ds\qq\q~~~+\Big\{C(s)X(s)+\bar C(s)\dbE[X(s)]+D(s)\big\{\Th^*(s)X(s)+\bar\Th^*(s)\dbE[X(s)]+v(s)\big\}\\
\ns\ds\qq\qq\qq~+\bar D(s)\dbE\big\{\Th^*(s)X(s)+\bar\Th^*(s)\dbE[X(s)]+v(s)\big\}+\si(s)\Big\}dW(s),\qq s\in[t,T],\\
\ns\ds X(t)=\xi.\ea\right.$$
Therefore,
$$\ba{ll}
\ds J(t,\xi;\Th^*(\cd)X^*(\cd)+\bar\Th^*(\cd)\dbE[X^*(\cd)]+v^*(\cd))\\
\ns\ds\les J(t,\xi;\Th^*(\cd)X(\cd)+\bar\Th^*(\cd)\dbE[X(\cd)]+v(\cd))=J(t,\xi;u(\cd)).
\ea$$

(iii) $\Ra$ (i). For any $\xi\in L^2_{\cF_t}(\Om;\dbR^n)$ and $(\Th(\cd),\bar\Th(\cd),v(\cd))\in\sC[t,T]$,
let $X(\cd)$ be the solution to the following MF-SDE:
$$\left\{\2n\ba{ll}
\ds dX(s)=\Big\{A(s)X(s)+\bar A(s)\dbE[X(s)]+B(s)\big\{\Th(s)X(s)+\bar\Th(s)\dbE[X(s)]+v(s)\big\}\\
\ns\ds\qq\qq\qq~+\bar B(s)\dbE\big\{\Th(s)X(s)+\bar\Th(s)\dbE[X(s)]+v(s)\big\}+b(s)\Big\}ds\\
\ns\ds\qq\q~~~
+\Big\{C(s)X(s)+\bar C(s)\dbE[X(s)]+D(s)\big\{\Th(s)X(s)+\bar\Th(s)\dbE[X(s)]+v(s)\big\}\\
\ns\ds\qq\qq\qq~+\bar D(s)\dbE\big\{\Th(s)X(s)+\bar\Th(s)\dbE[X(s)]+v(s)\big\}+\si(s)\Big\}dW(s),\qq s\in[t,T],\\
\ns\ds X(t)=\xi.\ea\right.$$
Set
$$u(\cd)\deq\Th(\cd)X(\cd)+\bar\Th(\cd)\dbE[X(\cd)]+v(\cd)\in\cU[t,T].$$
Then, by uniqueness, $X(\cd)$ also solves MF-SDE \rf{12.25-1}. Thus,
$$\ba{ll}
\ds J(t,\xi;\Th^*(\cd)X^*(\cd)+\bar\Th^*(\cd)\dbE[X^*(\cd)]+v^*(\cd))\\
\ns\ds\les J(t,\xi;u(\cd))=J(t,\xi;\Th(\cd)X(\cd)+\bar\Th(\cd)\dbE[X(\cd)]+v(\cd)).
\ea$$
This completes the proof. \endpf

\ms

With the same proof, we have the following result for the weakly optimal closed-loop strategies.

\bp{equivalence*}\sl Let {\rm(H1)--(H2)} hold and let $(\Th^*(\cd),\bar\Th^*(\cd),v^*(\cd))\in\sC[t,T]$.
Then the following statements are equivalent:
\begin{enumerate}[~~\,\rm(i)]
\item $(\Th^*(\cd),\bar\Th^*(\cd),v^*(\cd))$ is a weakly optimal closed-loop strategy
of Problem {\rm(MF-LQ)} on $[t,T]$;
\item The following holds:
$$\ba{ll}
\ds J(t,x;\Th^*(\cd)X^*(\cd)+\bar\Th^*(\cd)\dbE[X^*(\cd)]+v^*(\cd))
\les J(t,x;\Th^*(\cd)X(\cd)+\bar\Th^*(\cd)\dbE[X(\cd)]+v(\cd)),\\
\ns\ds\qq\qq\qq\qq\qq\qq\qq\qq\qq\qq\qq\qq\qq~\forall(x,v(\cd))\in\dbR^n\times\cU[t,T],\ea$$
where $X^*(\cd)$ and $X(\cd)$ are the closed-loop state processes corresponding to
$(\Th^*(\cd),\bar\Th^*(\cd),v^*(\cd),x)$ and $(\Th^*(\cd),\bar\Th^*(\cd),v(\cd),x)$,
respectively;
\item The following holds:
$$J(t,x;\Th^*(\cd)X^*(\cd)+\bar\Th^*(\cd)\dbE[X^*(\cd)]+v^*(\cd))
\les J(t,x;u(\cd)),\q\forall(x,u(\cd))\in\dbR^n\times\cU[t,T],$$
where $X^*(\cd)$ is the closed-loop state process corresponding to
$(\Th^*(\cd),\bar\Th^*(\cd),v^*(\cd))$ and $x\in\dbR^n$.
\end{enumerate}
\ep

\br{br-1}\rm (i) An optimal open-loop control is allowed to depend on the initial state,
whereas an optimal closed-loop strategy is required to be independent of the initial state.

\ms

(ii) It is clear from Proposition \ref{equivalence} (iii) that the outcome control
$u^*(\cd)\equiv\Th^*(\cd)X^*(\cd)+\bar\Th^*(\cd)\dbE[X^*(\cd)]+v^*(\cd)$ of an optimal
closed-loop strategy $(\Th^*(\cd),\bar\Th^*(\cd),v^*(\cd))$ is an optimal open-loop control
of Problem (MF-LQ) for the initial pair $(t,X^*(t))$. Hence, closed-loop solvability implies
open-loop solvability.

\ms


(iii) Obviously, an optimal closed-loop strategy on $[t,T]$ is also weakly optimal.
For the classical LQ optimal control problems where $\dbE[X(\cd)]$ and $\dbE[u(\cd)]$ are absent,
it can be shown, using the results from Sun--Yong \cite{Sun-Yong 2014} and a completion of squares
technique, that the two concepts coincide. But for Problem (MF-LQ), the existence of a weakly optimal
closed-loop strategy does {\it not} guarantee the existence of an optimal closed-loop strategy.
To see this, we present the following example.
\er

\bex{Ex1}\rm Consider the following one-dimensional state equation
\bel{Ex1-state}\left\{\2n\ba{ll}
\ds dX(s)=\big\{u(s)-\dbE[u(s)]\big\}ds+\dbE[u(s)]dW(s),\qq s\in[t,1],\\
\ns\ds X(t)=\xi,
\ea\right.\ee
and cost functional
$$J(t,\xi;u(\cd))=\dbE\big[X(1)^2\big]+\big(\dbE[X(1)]\big)^2.$$
For any $x\in\dbR$ and $u(\cd)\in\cU[t,1]$, we have
$$J(t,x;u(\cd))=\dbE\big[X(1)^2\big]+\big(\dbE[X(1)]\big)^2\ges 2\big(\dbE[X(1)]\big)^2=2x^2.$$
On the other hand, it is clear that $(\Th^*(\cd),\bar\Th^*(\cd),v^*(\cd))\equiv(0,0,0)$ satisfies
$$J(t,x;\Th^*(\cd)X^*(\cd)+\bar\Th^*(\cd)\dbE[X^*(\cd)]+v^*(\cd))=J(t,x;0)=2x^2,
\qq\forall x\in\dbR.$$
Thus, by Proposition \ref{equivalence*}, $(0,0,0)$ is a weakly optimal closed-loop strategy of
the problem on $[t,1]$.

\ms

Let us now show that the above problem does not admit an optimal closed-loop strategy
on any $[t,1]$ with $0<t<1$. Assume the contrary; i.e.,
let $(\Th^*(\cd),\bar\Th^*(\cd),v^*(\cd))\in\sC[t,1]$ satisfy \rf{strong**} on $[t,1]$
for some $0<t<1$. For any $x\in\dbR$, take $\xi=W(t)x$ and $u(s)\equiv{W(t)x\over t-1}$.
The corresponding solution of \rf{Ex1-state} is
$$\ba{lll}
\ds X(s)\4n&=\4n&\ds\xi+\int_t^s\big\{u(r)-\dbE[u(r)]\big\}dr+\int_t^s\dbE[u(r)]dW(r)\\
\ns\4n&=\4n&\ds W(t)x+{s-t\over t-1}W(t)x,\qq s\in[t,1].
\ea$$
Note that $X(1)=0$. Thus, \rf{strong**} implies that
$$\ba{lll}
\ds\dbE\big[X^*(1)^2\big]+\big(\dbE[X^*(1)]\big)^2
\4n&=\4n&\ds J(t,W(t)x;\Th^*(\cd)X^*(\cd)+\bar\Th^*(\cd)\dbE[X^*(\cd)]+v^*(\cd))\\
\ns\4n&\les\4n&\ds J\bigg(t,W(t)x;{W(t)x\over t-1}\bigg)=0,\qq\forall x\in\dbR,
\ea$$
where $X^*(\cd)$ is the solution to the following closed-loop system:
$$\left\{\2n\ba{ll}
\ds dX^*(s)=\Big\{\Th^*(s)\big(X^*(s)-\dbE[X^*(s)]\big)+v^*(s)-\dbE[v^*(s)]\Big\}ds\\
\ns\ds\qq\qq~~+\Big\{[\Th^*(s)+\bar\Th^*(s)]\dbE[X^*(s)]+\dbE[v^*(s)]\Big\}dW(s),
\qq s\in[t,1],\\
\ns\ds X^*(t)=W(t)x.\ea\right.$$
It follows that
$$X^*(1)=0,\qq\forall x\in\dbR.$$
Note that $\dbE[X^*(s)]\equiv0$. Then,
$$\ba{lll}
\ds 0=X^*(1)=e^{\int_t^1\Th^*(s)ds}W(t)x
+\int_t^1e^{\int_r^1\Th^*(s)ds}\big\{v^*(r)-\dbE[v^*(r)]\big\}dr\\
\ns\ds\qq\qq\q~~+\int_t^1e^{\int_r^1\Th^*(s)ds}\dbE[v^*(r)]dW(r),
\qq\forall x\in\dbR.
\ea$$
But this is impossible since it has to be true for all $x\in\dbR$.
\ex

We conclude this section by introducing the coupled generalized Riccati equations
(GREs, for short), whose regular solvability will turn out to be necessary and
sufficient for the closed-loop solvability of Problem (MF-LQ) in the next two
sections.

\ms

The GREs associated with Problem (MF-LQ) are two coupled nonlinear differential equations of the
following form (for simplicity of notation, we will usually suppress the time variable $s$ below):
\bel{Ric}\left\{\2n\ba{ll}
\ds\dot P+PA+A^\top P+C^\top PC+Q\\
\ns\ds\q-\,\big(PB+C^\top PD+S^\top\big)\big(R+D^\top PD\big)^\dag\big(B^\top P+D^\top PC+S\big)=0,
\qq\ae~s\in[t,T],\\
\ns\ds\dot\Pi+\Pi(A+\bar A)+(A+\bar A)^\top\Pi+Q+\bar Q+(C+\bar C)^\top P(C+\bar C)\\
\ns\ds\q-\,\big[\Pi(B+\bar B)+(C+\bar C)^\top P(D+\bar D)+(S+\bar S)^\top\big]
\big[R+\bar R+(D+\bar D)^\top P(D+\bar D)\big]^\dag\\
\ns\ds\qq~\cd\big[(B+\bar B)^\top\Pi+(D+\bar D)^\top P(C+\bar C)+(S+\bar S)\big]=0,
\qq\ae~s\in[t,T],\\
\ns\ds P(T)=G,\qq \Pi(T)=G+\bar G.\ea\right.\ee

\bde{reg-sol}\rm A solution $(P(\cd),\Pi(\cd))\in C([t,T];\dbS^n)\times C([t,T];\dbS^n)$
of \rf{Ric} is said to be {\it regular} if
\bel{positive}\Si\equiv R+D^\top PD\ges0,\qq
\bar\Si\equiv R+\bar R+(D+\bar D)^\top P(D+\bar D)\ges0,\ee
\bel{range}\left\{\2n\ba{ll}
\ds\cR\big(B^\top P+D^\top PC+S\big)\subseteq\cR(\Si),\\
\ns\ds\cR\big((B+\bar B)^\top\Pi+(D+\bar D)^\top P(C+\bar C)+(S+\bar S)\big)\subseteq\cR(\bar\Si),
\ea\right.\ee
and
\bel{L2}\left\{\2n\ba{ll}
\ds\Si^\dag\big(B^\top P+D^\top PC+S\big)\in L^2(t,T;\dbR^{m\times n}),\\
\ns\ds\bar\Si^\dag\big[(B+\bar B)^\top\Pi+(D+\bar D)^\top P(C+\bar C)+(S+\bar S)\big]
\in L^2(t,T;\dbR^{m\times n}).\ea\right.\ee
The GREs \rf{Ric} is said to be {\it regularly solvable} on $[t,T]$ if it admits a regular solution.
\ede

\section{Necessary conditions for closed-loop solvability}

In this section we will deduce necessary conditions for the closed-loop solvability of Problem (MF-LQ).
In particular, we shall establish the necessity of the regular solvability of GREs \rf{Ric} by a matrix
minimum principle.

\ms

Let $\Th^*(\cd),\bar\Th^*(\cd)\in L^2(t,T;\dbR^{m\times n})$
and consider the following state equation
$$\left\{\2n\ba{ll}
\ds dX(s)=\Big\{AX+\bar A\dbE[X]+B\big(\Th^*X+\bar\Th^*\dbE[X]+u\big)
+\bar B\dbE\big(\Th^*X+\bar\Th^*\dbE[X]+u\big)+b\Big\}ds\\
\ns\ds\qq\q~~~+\Big\{CX+\bar C\dbE[X]+D\big(\Th^* X+\bar\Th^*\dbE[X]+u\big)
+\bar D\dbE\big(\Th^*X+\bar\Th^*\dbE[X]+u\big)+\si\Big\}dW(s),\\
\ns\ds X(t)=\xi,\ea\right.$$
and cost functional
$$\ba{ll}
\ds\wt J(t,\xi;u(\cd))\deq J(t,\xi;\Th^*(\cd) X(\cd)+\bar\Th^*(\cd)\dbE[X(\cd)]+u(\cd))\\
\ns\ds=\dbE\Bigg\{\lan GX(T),X(T)\ran+2\lan g,X(T)\ran
+\blan\bar G\dbE[X(T)],\dbE[X(T)]\bran+2\lan\bar g,\dbE[X(T)]\ran\\
\ns\ds\qq~+\int_t^T\Bigg[\llan\begin{pmatrix}Q&S^\top\\S&R\end{pmatrix}
                              \begin{pmatrix}X\\ \Th^*X+\bar\Th^*\dbE[X]+u\end{pmatrix},
                              \begin{pmatrix}X\\ \Th^*X+\bar\Th^*\dbE[X]+u\end{pmatrix}\rran\\
\ns\ds\qq\qq\qq~+2\llan\begin{pmatrix}q\\ \rho\end{pmatrix},
                       \begin{pmatrix}X\\ \Th^*X+\bar\Th^*\dbE[X]+u\end{pmatrix}\rran\Bigg]ds\\
\ns\ds\qq~+\int_t^T\Bigg[\llan\begin{pmatrix}\bar Q&\bar S^\top\\\bar S&\bar R\end{pmatrix}
                              \begin{pmatrix}\dbE[X]\\(\Th^*+\bar\Th^*)\dbE[X]+\dbE[u]\end{pmatrix},
                              \begin{pmatrix}\dbE[X]\\(\Th^*+\bar\Th^*)\dbE[X]+\dbE[u]\end{pmatrix}\rran\\
\ns\ds\qq\qq\qq~+2\llan\begin{pmatrix}\bar q\\ \bar\rho\end{pmatrix},
                       \begin{pmatrix}\dbE[X]\\(\Th^*+\bar\Th^*)\dbE[X]+\dbE[u]\end{pmatrix}\rran\Bigg]ds\Bigg\}\\
\ns\ds=\dbE\Bigg\{\lan GX(T),X(T)\ran+2\lan g,X(T)\ran
+\blan\bar G\dbE[X(T)],\dbE[X(T)]\bran+2\lan\bar g,\dbE[X(T)]\ran\\
\ns\ds\qq~+\int_t^T\lt[\llan\begin{pmatrix}\wt Q&\wt S^\top\\ \wt S&R\end{pmatrix}
                            \begin{pmatrix}X\\ u\end{pmatrix},
                            \begin{pmatrix}X\\ u\end{pmatrix}\rran
+2\llan\begin{pmatrix}\wt q\\ \rho\end{pmatrix},
       \begin{pmatrix}X\\ u\end{pmatrix}\rran\rt]ds\\
\ns\ds\qq~+\int_t^T\lt[\llan\begin{pmatrix}\h Q&\h S^\top\\ \h S&\bar R\end{pmatrix}
                            \begin{pmatrix}\dbE[X]\\ \dbE[u]\end{pmatrix},
                            \begin{pmatrix}\dbE[X]\\ \dbE[u]\end{pmatrix}\rran
+2\llan\begin{pmatrix}\h q\\ \bar\rho\end{pmatrix},
       \begin{pmatrix}\dbE[X]\\\dbE[u]\end{pmatrix}\rran\rt]ds\Bigg\},\ea$$
where
$$\left\{\2n\ba{cll}
\wt Q=Q\4n&+\4n&(\Th^*)^\top S+S^\top\Th^*+(\Th^*)^\top R\Th^*,
\q \wt S=S+R\Th^*,\q\wt q=q+(\Th^*)^\top\rho,\\
\ns\h Q=\bar Q\4n&+\4n&(\Th^*+\bar\Th^*)^\top\bar S+\bar S^\top(\Th^*+\bar\Th^*)
+(\Th^*+\bar\Th^*)^\top\bar R(\Th^*+\bar\Th^*)\\
\ns\4n&+\4n&(\bar\Th^*)^\top R\bar\Th^*+(\bar\Th^*)^\top S+S^\top\bar\Th^*
+(\bar\Th^*)^\top R\Th^*+(\Th^*)^\top R\bar\Th^*,\\
\ns\h S=\bar S\4n&+\4n&\bar R(\Th^*+\bar \Th^*)+R\bar \Th^*,
\q\h q=\bar q+(\Th^*+\bar\Th^*)^\top\bar\rho+(\bar\Th^*)^\top\dbE[\rho].
\ea\right.$$
By Proposition \ref{equivalence} (ii), $(\Th^*(\cd),\bar\Th^*(\cd),u^*(\cd))\in\sC[t,T]$ is
an optimal closed-loop strategy of Problem (MF-LQ) on $[t,T]$ if and only if for any
$\xi\in L^2_{\cF_t}(\Om;\dbR^n)$, $u^*(\cd)$ is an optimal open-loop control of the problem
with the above state equation and cost functional. This leads to the following result.

\bp{prop2}\sl Let {\rm(H1)--(H2)} hold. If $(\Th^*(\cd),\bar\Th^*(\cd),u^*(\cd))\in\sC[t,T]$
is an optimal closed-loop strategy of Problem {\rm(MF-LQ)} on $[t,T]$,
then $(\Th^*(\cd),\bar\Th^*(\cd),0)$ is an optimal closed-loop strategy of Problem
{\rm(MF-LQ)}$^0$ on $[t,T]$.
\ep

\it Proof. \rm By the preceding discussion and \cite[Theorem 2.3]{Sun},
we see that $(\Th^*(\cd),\bar\Th^*(\cd),u^*(\cd))$ is an optimal closed-loop strategy of
Problem (MF-LQ) on $[t,T]$ if and only if for any $\xi\in L^2_{\cF_t}(\Om;\dbR^n)$,
the adapted solution $(X^*(\cd),Y^*(\cd),Z^*(\cd))$ to the following mean-field forward-backward
stochastic differential equation (MF-FBSDE, for short):
\bel{MF-FBSDE}\left\{\2n\ba{ll}
\ds dX^*(s)=\Big\{(A+B\Th^*)X^*+\big[\bar A+B\bar\Th^*+\bar B(\Th^*+\bar\Th^*)\big]\dbE[X^*]
+Bu^*+\bar B\dbE[u^*]+b\Big\}ds\\
\ns\ds\qq\qq~+\Big\{(C+D\Th^*)X^*+\big[\bar C+D\bar\Th^*+\bar D(\Th^*+\bar\Th^*)\big]\dbE[X^*]
+Du^*+\bar D\dbE[u^*]+\si\Big\}dW(s),\\
%
%
\ns\ds dY^*(s)=-\Big\{(A+B\Th^*)^\top Y^*+\big[\bar A+B\bar\Th^*+\bar B(\Th^*+\bar\Th^*)\big]^\top\dbE[Y^*]\\
\ns\ds\qq\qq\qq\2n~+(C+D\Th^*)^\top Z^*+\big[\bar C+D\bar\Th^*+\bar D(\Th^*+\bar\Th^*)\big]^\top\dbE[Z^*]\\
\ns\ds\qq\qq\qq\2n~+\wt QX^*+\h Q\dbE[X^*]+\wt S^\top u^*+\h S^\top\dbE[u^*]+\wt q+\h q\Big\}ds+Z^*dW(s),
\qq s\in[t,T],\\
\ns\ds X^*(t)=\xi,\qq Y^*(T)=GX^*(T)+\bar G\dbE[X^*(T)]+g+\bar g,
\ea\right.\ee
satisfies
\bel{yueshu}Ru^*+ B^\top Y^*+ D^\top Z^*+\wt SX^*+\rho
+\bar R\dbE[u^*]+\bar B^\top\dbE[Y^*]+\bar D^\top\dbE[Z^*]+\h S\dbE[X^*]+\bar\rho=0,\ee
and the following condition hold:
$$\ba{ll}
\ds\dbE\Bigg\{\lan GX(T),X(T)\ran+\blan\bar G\dbE[X(T)],\dbE[X(T)]\bran
+\int_t^T\[\blan\wt QX,X\bran+2\blan\wt SX,u\bran+\lan Ru,u\ran\]ds\\
\ns\ds\q~+\int_t^T\[\blan\h Q\dbE[X],\dbE[X]\bran+2\blan\h S\dbE[X],\dbE[u]\bran
+\blan\bar R\dbE[u],\dbE[u]\bran\]ds\Bigg\}\ges0,\qq\forall u(\cd)\in\cU[t,T],
\ea$$
where $X(\cd)$ is the solution of
$$\left\{\2n\ba{ll}
\ds dX(s)=\Big\{(A+B\Th^*)X+\big[\bar A+B\bar\Th^*\1n+\bar B(\Th^*\1n+\bar\Th^*)\big]\dbE[X]
+Bu+\bar B\dbE[u]\Big\}ds \\
\ns\ds\qq\qq~+\Big\{(C+D\Th^*)X+\big[\bar C+D\bar\Th^*\1n+\bar D(\Th^*\1n+\bar\Th^*)\big]\dbE[X]
+Du+\bar D\dbE[u]\Big\}dW(s),\qq s\in[t,T], \\
\ns\ds X(t)=0.\ea\right.$$
Since the MF-FBSDE \rf{MF-FBSDE} admits a solution for each $\xi\in L^2_{\cF_t}(\Om;\dbR^n)$
and $(\Th^*(\cd),\bar\Th^*(\cd),u^*(\cd))$ is independent of $\xi$, by subtracting solutions
corresponding $\xi$ and $0$, the later from the former, we see that for any
$\xi\in L^2_{\cF_t}(\Om;\dbR^n)$, the following MF-FBSDE:
$$\left\{\2n\ba{ll}
\ds dX(s)=\Big\{(A+B\Th^*)X+\big[\bar A+B\bar\Th^*+\bar B(\Th^*+\bar\Th^*)\big]\dbE[X]\Big\}ds\\
\ns\ds\qq\q~~~
+\Big\{(C+D\Th^*)X+\big[\bar C+D\bar\Th^*+\bar D(\Th^*+\bar\Th^*)\big]\dbE[X]\Big\}dW(s),\qq s\in[t,T],\\
\ns\ds
dY(s)=-\Big\{(A+B\Th^*)^\top Y+\big[\bar A+B\bar\Th^*+\bar B(\Th^*+\bar\Th^*)\big]^\top\dbE[Y]+(C+D\Th^*)^\top Z\\
\ns\ds\qq\qq\q~
+\big[\bar C+D\bar\Th^*+\bar D(\Th^*+\bar\Th^*)\big]^\top\dbE[Z]+\wt QX+\h Q\dbE[X]\Big\}ds+ZdW(s),\qq s\in[t,T],\\
\ns\ds X(t)=\xi,\qq Y(T)=GX(T)+\bar G\dbE[X(T)],
\ea\right.$$
also admits an adapted solution $(X(\cd),Y(\cd),Z(\cd))$ satisfying
$$B^\top Y+ D^\top Z+\wt SX+\bar B^\top\dbE[Y]+\bar D^\top\dbE[Z]+\h S\dbE[X]=0.$$
It follows, again from \cite[Theorem 2.3]{Sun}, that $(\Th^*(\cd),\bar\Th^*(\cd),0)$ is an
optimal closed-loop strategy of Problem {\rm(MF-LQ)}$^0$ on $[t,T]$.
\endpf

\ms

Now let us look at Problem {\rm(MF-LQ)}$^0$. If we consider only closed-loop strategies
of the form $(\Th(\cd),\bar\Th(\cd),0)$, then the state equation becomes
$$\left\{\2n\ba{ll}
\ds dX(s)=\Big\{(A+B\Th)X+\big[\bar A+\bar B\Th+(B+\bar B)\bar\Th\big]\dbE[X]\Big\}ds\\
\ns\ds\qq\q~~~+\Big\{(C+D\Th)X+\big[\bar C+\bar D\Th+(D+\bar D)\bar\Th\big]\dbE[X]\Big\}dW(s),
\qq s\in[t,T], \\
\ns\ds X(t)=\xi,\ea\right.$$
and $\dbE[X(\cd)]$ satisfies
$$\left\{\2n\ba{ll}
\ds d\dbE[X(s)]=\big[A+\bar A+(B+\bar B)(\Th+\bar\Th)\big]\dbE[X]ds,
\qq s\in[t,T], \\
\ns\ds \dbE[X(t)]=\dbE[\xi].\ea\right.$$
By It\^o's formula, the matrices $\BX(s)\deq\dbE\big[X(s)X(s)^\top\big]$ and
$\BY(s)\deq\dbE[X(s)]\dbE[X(s)]^\top$ satisfy the matrix-valued ordinary differential
equations (ODEs, for short)
\bel{BX}\left\{\2n\ba{ll}
\ds \dot{\BX}=(A+B\Th)\BX+\BX(A+B\Th)^\top\1n+(C+D\Th)\BX(C+D\Th)^\top\\
\ns\ds\qq~+\big[\bar A+\bar B\Th+(B+\bar B)\bar\Th\big]\BY
+\BY\big[\bar A+\bar B\Th+(B+\bar B)\bar\Th\big]^\top\\
\ns\ds\qq~+(C+D\Th)\BY\big[\bar C+\bar D\Th+(D+\bar D)\bar\Th\big]^\top\1n
+\big[\bar C+\bar D\Th+(D+\bar D)\bar\Th\big]\BY(C+D\Th)^\top\\
\ns\ds\qq~+\big[\bar C+\bar D\Th+(D+\bar D)\bar\Th\big]\BY
\big[\bar C+\bar D\Th+(D+\bar D)\bar\Th\big]^\top, \qq s\in[t,T],\\
\ns\ds \BX(t)=\dbE[\xi\xi^\top],\ea\right.\ee
and
\bel{BY}\left\{\2n\ba{ll}
\ds \dot{\BY}=\big[A+\bar A+(B+\bar B)(\Th+\bar\Th)\big]\BY
+\BY\big[A+\bar A+(B+\bar B)(\Th+\bar\Th)\big]^\top, \qq s\in[t,T],\\
\ns\ds \BY(t)=\dbE[\xi]\dbE[\xi]^\top,\ea\right.\ee
respectively. The cost functional $J^0(t,\xi;\Th(\cd)X(\cd)+\bar\Th(\cd)\dbE[X(\cd)])$
can be expressed equivalently as
\bel{9-25-cost}\BJ(t,\xi;\Th(\cd),\bar\Th(\cd))=\tr\big[G\BX(T)+\bar G\BY(T)\big]
+\int_t^T\tr\big[M(s)\BX(s)+N(s)\BY(s)\big]ds,\ee
where
$$\left\{\2n\ba{cll}
M=Q\4n&+\4n&\Th^\top S+S^\top\Th+\Th^\top R\Th, \\
\ns N=\bar Q\4n&+\4n&(\Th+\bar\Th)^\top\bar S+\bar S^\top(\Th+\bar\Th)+(\Th+\bar\Th)^\top\bar R(\Th+\bar\Th)\\
\ns\4n&+\4n&\bar\Th^\top R\bar\Th+\bar\Th^\top S+S^\top\bar\Th+\bar\Th^\top R\Th+\Th^\top R\bar\Th.
\ea\right.$$
Then we may pose the following deterministic optimal control problem.

\ms

\bf Problem (O). \rm For any given $(t,\xi)\in[0,T)\times L^2_{\cF_t}(\Om;\dbR^n)$,
find $\Th^*(\cd),\bar\Th^*(\cd)\in L^2(t,T;\dbR^{m\times n})$ such that
$$\ba{ll}
\ds \BJ(t,\xi;\Th^*(\cd),\bar\Th^*(\cd))\les\BJ(t,\xi;\Th(\cd),\bar\Th(\cd)),
\qq\forall \Th(\cd),\bar\Th(\cd)\in L^2(t,T;\dbR^{m\times n}).\ea$$

\ms

Rewrite \rf{BX}--\rf{BY} as
$$\left\{\2n\ba{ll}
\ds\begin{pmatrix}\dot\BX(s)\\ \dot\BY(s)\end{pmatrix}
=\begin{pmatrix}F_1(\BX(s),\BY(s),\Th(s),\bar\Th(s),s)\\
                F_2(\BY(s),\Th(s),\bar\Th(s),s)\end{pmatrix},\qq s\in[t,T],\\
\ns\ds\BX(t)=\dbE[\xi\xi^\top],\q \BY(t)=\dbE[\xi]\dbE[\xi]^\top,
\ea\right.$$
and denote the integrand in \rf{9-25-cost} by $L(\BX(s),\BY(s),\Th(s),\bar\Th(s),s)$.
We present the following matrix minimum principle for Problem (O).
The interested reader is referred to Athans \cite{Athans 1968} for a proof.

\bl{lmm-MMP}\sl Let {\rm(H1)--(H2)} hold. Suppose that $(\Th^*(\cd),\bar\Th^*(\cd))$
is an optimal control of Problem {\rm(O)} for the initial pair $(t,\xi)$
and let $(\BX^*(\cd),\BY^*(\cd))$ be the corresponding optimal state process.
Then there exist matrix-valued functions $P(\cd)$ and $\L(\cd)$ satisfying the following
ODEs (the variable $s\in[t,T]$ is suppressed)
\bel{MMP-equ}\left\{\2n\ba{ll}
\ds\begin{pmatrix}\dot P\\ \dot\L\end{pmatrix}=-\begin{pmatrix}
{\pr\over\pr \BX^*}L(\BX^*,\BY^*,\Th^*,\bar\Th^*)
+{\pr\over\pr \BX^*}\tr\big[F_1(\BX^*,\BY^*,\Th^*,\bar\Th^*)P^\top+F_2(\BY^*,\Th^*,\bar\Th^*)\L^\top\big]\\
{\pr\over\pr \BY^*}L(\BX^*,\BY^*,\Th^*,\bar\Th^*)
+{\pr\over\pr \BY^*}\tr\big[F_1(\BX^*,\BY^*,\Th^*,\bar\Th^*)P^\top+F_2(\BY^*,\Th^*,\bar\Th^*)\L^\top\big]
\end{pmatrix},\\
\ns\ds P(T)=G,\qq \L(T)=\bar G,\ea\right.\ee
with constraints
\bel{MMP-cons}\left\{\2n\ba{ll}
\no\ss\ds{\pr\over\pr \Th^*}L(\BX^*,\BY^*,\Th^*,\bar\Th^*)
+{\pr\over\pr \Th^*}\tr\big[F_1(\BX^*,\BY^*,\Th^*,\bar\Th^*)P^\top+F_2(\BY^*,\Th^*,\bar\Th^*)\L^\top\big]=0,\\
\ns\ds {\pr\over\pr \bar\Th^*}L(\BX^*,\BY^*,\Th^*,\bar\Th^*)
+{\pr\over\pr \bar\Th^*}\tr\big[F_1(\BX^*,\BY^*,\Th^*,\bar\Th^*)P^\top+F_2(\BY^*,\Th^*,\bar\Th^*)\L^\top\big]=0.
\ea\right.\ee
\el

Now, we are ready to state and prove the principal result of this section.

\bt{Th-Nec}\sl Let {\rm(H1)--(H2)} hold and $t\in(0,T)$. If Problem {\rm(MF-LQ)} admits an
optimal closed-loop strategy on $[t,T]$, then the GREs \rf{Ric} is regularly solvable on $[t,T]$.
\et

\it Proof. \rm Suppose that $(\Th^*(\cd),\bar\Th^*(\cd),u^*(\cd))\in\sC[t,T]$ is an optimal closed-loop
strategy of Problem (MF-LQ) on $[t,T]$. Then, by Proposition \ref{prop2},
$(\Th^*(\cd),\bar\Th^*(\cd),0)$ is an optimal closed-loop strategy of Problem (MF-LQ)$^0$
on $[t,T]$, and it follows from Definition \ref{optimal closed-loop strategy} (i) that $(\Th^*(\cd),\bar\Th^*(\cd))$
is an optimal control of Problem (O) for any $\xi\in L^2_{\cF_t}(\Om;\dbR^n)$.
Thus, by the matrix minimum principle, Lemma \ref{lmm-MMP}, there exist functions
$P(\cd),\L(\cd):[t,T]\to\dbR^n$ such that \rf{MMP-equ}--\rf{MMP-cons} hold.
By a straightforward calculation, we see from the first equation in \rf{MMP-equ} that
$P(\cd)$ satisfies
\bel{Ri-P}\left\{\2n\ba{ll}
\ds \dot P+(A+B\Th^*)^\top P+P(A+B\Th^*)+(C+D\Th^*)^\top P(C+D\Th^*)\\
\ns\ds\q+\,Q+(\Th^*)^\top S+S^\top\Th^*+(\Th^*)^\top R\Th^*=0,\\
\ns\ds P(T)=G,\ea\right.\ee
and from the second equation in \rf{MMP-equ}, we see that $\L(\cd)$ satisfies
\bel{Ri-L}\left\{\2n\ba{ll}
\ds 0=\dot\L+\bar Q+(\Th^*+\bar\Th^*)^\top\bar S+\bar S^\top(\Th^*+\bar\Th^*)
+(\Th^*+\bar\Th^*)^\top\bar R(\Th^*+\bar\Th^*)\\
\ns\ds\qq~~\1n+(\bar\Th^*)^\top R\bar\Th^*+(\bar\Th^*)^\top S+S^\top\bar\Th^*
+(\bar\Th^*)^\top R\Th^*+(\Th^*)^\top R\bar\Th^*\\
\ns\ds\qq~~\1n+\big[\bar A+\bar B\Th^*+(B+\bar B)\bar\Th^*\big]^\top P
+P\big[\bar A+\bar B\Th^*+(B+\bar B)\bar\Th^*\big]\\
\ns\ds\qq~~\1n+\big[C+D\Th^*\big]^\top P\big[\bar C+\bar D\Th^*+(D+\bar D)\bar\Th^*\big]\\
\ns\ds\qq~~\1n+\big[\bar C+\bar D\Th^*+(D+\bar D)\bar\Th^*\big]^\top P\big[C+D\Th^*\big]\\
\ns\ds\qq~~\1n+\big[\bar C+\bar D\Th^*+(D+\bar D)\bar\Th^*\big]^\top P
\big[\bar C+\bar D\Th^*+(D+\bar D)\bar\Th^*\big]\\
\ns\ds\qq~~\1n+\big[A+\bar A+(B+\bar B)(\Th^*+\bar\Th^*)\big]^\top\L
+\L\big[A+\bar A+(B+\bar B)(\Th^*+\bar\Th^*)\big],\\
\ns\ds \L(T)=\bar G.
\ea\right.\ee
Note that $P(\cd)^\top$ and $\L(\cd)^\top$ also solve \rf{Ri-P} and \rf{Ri-L}, respectively.
Hence, by uniqueness, we have $P(\cd)=P(\cd)^\top$ and $\L(\cd)=\L(\cd)^\top$. Let
$$\Pi(\cd)=P(\cd)+\L(\cd),\qq \D(\cd)=\Th^*(\cd)+\bar\Th^*(\cd).$$
Then, $\Pi(T)=G+\bar G$ and
\bel{Ri-Pi}\ba{lllll}
\ds 0\4n&=\4n&\dot\Pi\5n&+\5n&Q+\bar Q+\D^\top(S+\bar S)+(S+\bar S)^\top\D+\D^\top(R+\bar R)\D\\
\ns\4n&~\4n&\5n&+\5n&\big[A+\bar A+(B+\bar B)\D\big]^\top P+P\big[A+\bar A+(B+\bar B)\D\big]\\
\ns\4n&~\4n&\5n&+\5n&\big[C+D\Th^*\big]^\top P\big[C+\bar C+(D+\bar D)\D\big]\\
\ns\4n&~\4n&\5n&+\5n&\big[\bar C+\bar D\Th^*+(D+\bar D)\bar\Th^*\big]^\top P\big[C+\bar C+(D+\bar D)\D\big]\\
\ns\4n&~\4n&\5n&+\5n&\big[A+\bar A+(B+\bar B)\D\big]^\top\L+\L\big[A+\bar A+(B+\bar B)\D\big]\\
\ns\4n&=\4n&\dot\Pi\5n&+\5n&Q+\bar Q+\D^\top(S+\bar S)+(S+\bar S)^\top\D+\D^\top(R+\bar R)\D\\
\ns\4n&~\4n&\5n&+\5n&\big[A+\bar A+(B+\bar B)\D\big]^\top\Pi+\Pi\big[A+\bar A+(B+\bar B)\D\big]\\
\ns\4n&~\4n&\5n&+\5n&\big[C+\bar C+(D+\bar D)\D\big]^\top P\big[C+\bar C+(D+\bar D)\D\big]\\
\ns\4n&=\4n&\dot\Pi\5n&+\5n&\big[A+\bar A+(B+\bar B)\D\big]^\top\Pi+\Pi\big[A+\bar A+(B+\bar B)\D\big]\\
\ns\4n&~\4n&\5n&+\5n&Q+\bar Q+(C+\bar C)^\top P(C+\bar C)+\D^\top\big[R+\bar R+(D+\bar D)^\top P(D+\bar D)\big]\D\\
\ns\4n&~\4n&\5n&+\5n&\D^\top\big[(D+\bar D)^\top P(C+\bar C)+S+\bar S\big]
+\big[(D+\bar D)^\top P(C+\bar C)+S+\bar S\big]^\top\D.
\ea\ee
Also, from the first equality in \rf{MMP-cons}, we have (noting that $\BX^*$ and $\BY^*$ are symmetric)
\bel{cons-P}\ba{lll}
\ds 0\4n&=&\4n\ds 2S\BX^*+2R\Th^*\BX^*+2\bar S\BY^*+2\bar R\Th^*\BY^*
+2\bar R\bar\Th^*\BY^*+2R\bar\Th^*\BY^*\\
\ns\4n&~&\4n\ds+\,2B^\top P\BX^*+2D^\top PC\BX^*+2D^\top PD\Th^*\BX^*\\
\ns\4n&~&\4n\ds+\,2\bar B^\top P\BY^*+2\bar D^\top PC\BY^*
+2D^\top P\big[\bar C+(D+\bar D)\bar\Th^*\big]\BY^*\\
\ns\4n&~&\4n\ds+\,2D^\top P\bar D\Th^*\BY^*+2\bar D^\top PD\Th^*\BY^*
+2\bar D^\top P\big[\bar C+(D+\bar D)\bar\Th^*\big]\BY^*\\
\ns\4n&~&\4n\ds+\,2\bar D^\top P\bar D\Th^*\BY^*+2(B+\bar B)^\top\L \BY^*\\
\ns\4n&=&\4n\ds2\big[(R+D^\top PD)\Th^*+B^\top P+D^\top PC+S\big]\BX^*\\
\ns\4n&~&\4n\ds+\,2\Big\{(R+\bar R)\bar\Th^*+(\bar R+\bar D^\top PD)\Th^*
+\bar B^\top P+\bar D^\top PC+\bar S\\
\ns\4n&~&\4n\ds\qq~+(D+\bar D)^\top P\big[\bar C+\bar D\Th^*+(D+\bar D)\bar\Th^*\big]
+(B+\bar B)^\top\L \Big\}\BY^*\\
\ns\4n&=&\4n\ds2\big[(R+D^\top PD)\Th^*+B^\top P+D^\top PC+S\big]\BX^*\\
\ns\4n&~&\4n\ds+\,2\Big\{\1n-\big[(R+D^\top PD)\Th^*+B^\top P+D^\top PC+S\big]\\
\ns\4n&~&\4n\ds\qq\,+\,(R+\bar R)\D+(B+\bar B)^\top\Pi
+(D+\bar D)^\top P\big[C+\bar C+(D+\bar D)\D\big]+S+\bar S\Big\}\BY^*\\
\ns\4n&=&\4n\ds2\big[(R+D^\top PD)\Th^*+B^\top P+D^\top PC+S\big](\BX^*-\BY^*)\\
\ns\4n&~&\4n\ds+\,2\Big\{\big[R+\bar R+(D+\bar D)^\top P(D+\bar D)\big]\D+(B+\bar B)^\top\Pi\\
\ns\4n&~&\4n\ds\qq\,+\,(D+\bar D)^\top P(C+\bar C)+S+\bar S\Big\}\BY^*.
\ea\ee
Likewise, from the second equality in \rf{MMP-cons}, we have
\bel{cons-L}\ba{ll}
\ds 2\Big\{\big[R+\bar R+(D+\bar D)^\top P(D+\bar D)\big]\D+(B+\bar B)^\top\Pi\\
\ns\ds\q~+(D+\bar D)^\top P(C+\bar C)+S+\bar S\Big\}\BY^*=0.\ea\ee
Let $\F(\cd)$ be the solution to the $\dbR^{n\times n}$-valued ODE
$$\left\{\2n\ba{ll}
\ds\dot \F(s)=\wt A(s)\F(s),\qq s\in[t,T],\\
\ns\ds \F(t)=I,\ea\right.$$
where
$$\wt A\deq A+\bar A+(B+\bar B)(\Th^*+\bar\Th^*)=A+\bar A+(B+\bar B)\D.$$
Then
$$\BY^*(s)=\F(s)\dbE[\xi]\dbE[\xi]^\top\F(s)^\top,\qq s\in[t,T].$$
Denoting $\bar\Si\equiv R+\bar R+(D+\bar D)^\top P(D+\bar D)$, since \rf{cons-L} holds for all
$\xi\in L_{\cF_t}^2(\Om;\dbR^n)$ and $\F(s)$ is invertible for all $s\in[t,T]$, we must have
\bel{Pi-range}\bar\Si\D+(B+\bar B)^\top\Pi+(D+\bar D)^\top P(C+\bar C)+S+\bar S=0.\ee
Now take $\eta\in L_{\cF_t}^2(\Om;\dbR)$ with $\dbE\eta=0$ and $\dbE\eta^2=1$. Then for any $x\in\dbR^n$,
the trajectory $(X^*,Y^*)$ with respect to the initial state $\eta x$ satisfies $Y^*\equiv0$ and
$$\left\{\2n\ba{ll}
\ds dX^*(s)=\big[A(s)+B(s)\Th^*(s)\big]X^*(s)ds
+\big[C(s)+D(s)\Th^*(s)\big]X^*(s)dW(s),\qq s\in[t,T],\\
\ns\ds X^*(t)=\eta x.\ea\right.$$
Let $\Psi(\cd)$ be the solution to the following SDE for
$\dbR^{n\times n}$-valued process:
$$\left\{\2n\ba{ll}
\ds d\Psi(s)=\big[A(s)+B(s)\Th^*(s)\big]\Psi(s)ds
+\big[C(s)+D(s)\Th^*(s)\big]\Psi(s)dW(s),\qq s\in[t,T],\\
\ns\ds \Psi(t)=I.\ea\right.$$
Since $\Psi(s)$ is independent of $\cF_t$ for $s\in[t,T]$, we have
$$\BX^*(s)=\dbE\big[\Psi(s)x\eta^2x^\top\Psi(s)^\top\big]
=\dbE\big[\Psi(s)xx^\top\Psi(s)^\top\big],\qq s\in[t,T].$$
Hence, denoting $\Si\equiv R+D^\top PD$, we obtain from \rf{cons-P} that
$$\big(\Si\Th^*+B^\top P+D^\top PC+S\big)\dbE\big[\Psi xx^\top\Psi^\top\big]=0,
\qq\forall x\in\dbR^n,$$
which implies $\Si\Th^*+B^\top P+D^\top PC+S=0$.
It follows that $\cR(B^\top P+D^\top PC+S)\subseteq\cR(\Si)$.
Moreover, since $\Si^\dag\Si$ is an orthogonal projection, we have
$$\Si^\dag(B^\top P+D^\top PC+S)\in L^2(t,T;\dbR^{m\times n}),$$
and
\bel{Th-star}\Th^*=-\Si^\dag(B^\top P+D^\top PC+S)+(I-\Si^\dag\Si)\th,\ee
for some $\th(\cd)\in L^2(t,T;\dbR^{m\times n})$.
Similarly, from \rf{Pi-range} we have
$$\cR\big((B+\bar B)^\top\Pi+(D+\bar D)^\top P(C+\bar C)+(S+\bar S)\big)
\subseteq\cR(\bar \Si),$$
$$\bar\Si^\dag\big[(B+\bar B)^\top\Pi+(D+\bar D)^\top P(C+\bar C)+(S+\bar S)\big]
\in L^2(t,T;\dbR^{m\times n}),$$
and
\bel{Delta}\D=-\bar\Si^\dag\big[(B+\bar B)^\top\Pi+(D+\bar D)^\top P(C+\bar C)+(S+\bar S)\big]
+\big(I-\bar\Si^\dag\bar\Si\big)\t,\ee
for some $\t(\cd)\in L^2(t,T;\dbR^{m\times n})$. Substituting \rf{Th-star} and \rf{Delta}
back into \rf{Ri-P} and \rf{Ri-Pi}, respectively, we see that $(P(\cd),\Pi(\cd))$ satisfies the GREs
\rf{Ric}. In order to show that $(P(\cd),\Pi(\cd))$ is regular, it remains to prove that
$$\Si\equiv R+D^\top PD\ges0,\qq\bar\Si\equiv R+\bar R+(D+\bar D)^\top P(D+\bar D)\ges0.$$
For this we take any $u(\cd)\in\cU[t,T]$ and let $X(\cd)$ be the solution to
\bel{12.30-X}\left\{\2n\ba{ll}
\ds dX(s)=\Big\{A(s)X(s)+\bar A(s)\dbE[X(s)]+B(s)u(s)+\bar B(s)\dbE[u(s)]\Big\}ds\\
\ns\ds\qq\qq~
+\Big\{C(s)X(s)+\bar C(s)\dbE[X(s)]+D(s)u(s)+\bar D(s)\dbE[u(s)]\Big\}dW(s),\qq s\in[t,T],\\
\ns\ds X(t)=0.\ea\right.\ee
Applying It\^o's formula to $s\mapsto\lan P(s)(X(s)-\dbE[X(s)]),X(s)-\dbE[X(s)]\ran$
and $s\mapsto\lan\Pi(s)\dbE[X(s)],\dbE[X(s)]\ran$, we have
$$\ba{ll}
\ds J^0(t,0;u(\cd))\\
\ns\ds=
\dbE\Bigg\{\lan G(X(T)-\dbE[X(T)]),X(T)-\dbE[X(T)]\ran+\blan(G+\bar G)\dbE[X(T)],\dbE[X(T)]\bran\\
\ns\ds\qq~~
+\int_t^T\llan\begin{pmatrix}Q&S^\top \\ S&R\end{pmatrix}
              \begin{pmatrix}X-\dbE[X]\\ u-\dbE[u]\end{pmatrix},
              \begin{pmatrix}X-\dbE[X]\\ u-\dbE[u]\end{pmatrix}\rran ds\\
\ns\ds\qq~~
+\int_t^T\llan\begin{pmatrix}Q+\bar Q&(S+\bar S)^\top\\S+\bar S&R+\bar R\end{pmatrix}
              \begin{pmatrix}\dbE[X]\\ \dbE[u]\end{pmatrix},
              \begin{pmatrix}\dbE[X]\\ \dbE[u]\end{pmatrix}\rran ds\Bigg\}\\
\ns\ds=
\dbE\int_t^T\Big\{\blan\dot P(X-\dbE[X]),X-\dbE[X]\bran
+\blan P\big\{A(X-\dbE[X])+B(u-\dbE[u])\big\},X-\dbE[X]\bran\\
\ns\ds\qq\qq~
+\blan P(X-\dbE[X]),A(X-\dbE[X])+B(u-\dbE[u])\bran\\
\ns\ds\qq\qq~
+\blan P\big\{C(X-\dbE[X])+D(u-\dbE[u])+(C+\bar C)\dbE[X]+(D+\bar D)\dbE[u]\big\},\\
\ns\ds\qq\qq\qq\qq~
C(X-\dbE[X])+D(u-\dbE[u])+(C+\bar C)\dbE[X]+(D+\bar D)\dbE[u]\bran\Big\}ds\\
\ns\ds\q
+\,\int_t^T\Big\{\blan\dot\Pi\dbE[X],\dbE[X]\bran
+\blan\Pi\big\{(A+\bar A)\dbE[X]+(B+\bar B)\dbE[u]\big\},\dbE[X]\bran\\
\ns\ds\qq\qq~
+\blan\Pi\dbE[X],(A+\bar A)\dbE[X]+(B+\bar B)\dbE[u]\bran\Big\}ds\\
\ns\ds\q
+\,\dbE\int_t^T\Big\{\lan Q(X-\dbE[X]),X-\dbE[X]\ran
+2\lan S(X-\dbE[X]),u-\dbE[u]\ran+\lan R(u-\dbE[u]),u-\dbE[u]\ran\Big\}ds\\
\ns\ds\q
+\,\int_t^T\Big\{\blan(Q+\bar Q)\dbE[X],\dbE[X]\bran+2\blan(S+\bar S)\dbE[X],\dbE[u]\bran
+\blan(R+\bar R)\dbE[u],\dbE[u]\bran\Big\}ds\\
\ns\ds=\dbE\int_t^T\Big\{\blan\big(\dot P+PA+A^\top P+C^\top PC+Q\big)(X-\dbE[X]),X-\dbE[X]\bran\\
\ns\ds\qq\qq~
+2\blan\big(B^\top P+D^\top PC+S\big)(X-\dbE[X]),u-\dbE[u]\bran
+\blan\big(R+D^\top PD\big)(u-\dbE[u]),u-\dbE[u]\bran\Big\}ds\\
\ns\ds\q
+\,\int_t^T\Big\{\blan\big[\dot\Pi+\Pi(A+\bar A)+(A+\bar A)^\top\Pi+(C+\bar C)^\top P(C+\bar C)
+Q+\bar Q\big]\dbE[X],\dbE[X]\bran\\
\ns\ds\qq\qq~
+2\blan\big[(B+\bar B)^\top\Pi+(D+\bar D)^\top P(C+\bar C)+S+\bar S\,\big]\dbE[X],\dbE[u]\bran\\
\ns\ds\qq\qq~
+\blan\big[R+\bar R+(D+\bar D)^\top P(D+\bar D)\big]\dbE[u],\dbE[u]\bran\Big\}ds\\
\ns\ds=\dbE\int_t^T\Big\{\blan(\Th^*)^\top\Si\Th^*(X\1n-\1n\dbE[X]),X\1n-\1n\dbE[X]\bran
-2\lan\Si\Th^*(X\1n-\1n\dbE[X]),u\1n-\1n\dbE[u]\ran+\lan\Si(u\1n-\1n\dbE[u]),u\1n-\1n\dbE[u]\ran\Big\}ds\\
\ns\ds\q
+\,\int_t^T\Big\{\blan\D^\top\bar\Si\D\dbE[X],\dbE[X]\bran-2\blan\bar\Si\D\dbE[X],\dbE[u]\bran
+\blan\bar\Si\dbE[u],\dbE[u]\bran\Big\}ds\\
\ns\ds=\dbE\int_t^T\blan\Si\big\{u-\dbE[u]-\Th^*(X-\dbE[X])\big\},u-\dbE[u]-\Th^*(X-\dbE[X])\bran ds\\
\ns\ds\q+\,\int_t^T\blan\bar\Si\big(\dbE[u]-\D\dbE[X]\big),\dbE[u]-\D\dbE[X]\bran ds.
\ea$$
Since $(\Th^*(\cd),\bar\Th^*(\cd),0)$ is an optimal closed-loop strategy of
Problem (MF-LQ)$^0$ on $[t,T]$, we have
\bel{12.30}\ba{ll}
\ds \dbE\int_t^T\blan\Si\big\{u-\dbE[u]-\Th^*(X-\dbE[X])\big\},u-\dbE[u]-\Th^*(X-\dbE[X])\bran ds\\
\ns\ds\q~+\int_t^T\blan\bar\Si\big(\dbE[u]-\D\dbE[X]\big),\dbE[u]-\D\dbE[X]\bran ds\\
\ns\ds=J^0(t,0;u(\cd))\ges J^0(t,0;\Th^*(\cd)X^*(\cd)+\bar\Th^*(\cd)\dbE[X^*(\cd)])=0,
\qq\forall u(\cd)\in\cU[t,T].\ea\ee
Note that for any $u(\cd)\in\cU[t,T]$ of the form
$$u(s)=\Th^*(s)X(s)+v(s)W(s),\qq v(\cd)\in L^2(t,T;\dbR^m),$$
the corresponding solution $X(\cd)$ of \rf{12.30-X} satisfies $\dbE[X(\cd)]=0$
and hence $\dbE[u(\cd)]=0$. Then \rf{12.30} yields
$$\ba{lll}
0\4n&\les\4n&\ds\dbE\int_t^T\blan\Si(s)[u(s)-\Th^*(s)X(s)],u(s)-\Th^*(s)X(s)\bran ds\\
\ns\4n&=\4n&\ds \dbE\int_t^T\lan\Si(s)v(s)W(s),v(s)W(s)\ran ds\\
\ns\4n&\les\4n&\ds T\int_t^T\lan\Si(s)v(s),v(s)\ran ds,
\qq\forall v(\cd)\in L^2(t,T;\dbR^m),\ea$$
which implies that $\Si\ges0$. Likewise, for any $u(\cd)\in\cU[t,T]$ of the form
$$u(s)=\Th^*(s)\big\{X(s)-\dbE[X(s)]\big\}+\D(s)\dbE[X(s)]+v(s),\qq v(\cd)\in L^2(t,T;\dbR^m),$$
the corresponding solution $X(\cd)$ of \rf{12.30-X} satisfies
$$u(s)-\dbE[u(s)]=\Th^*(s)\big\{X(s)-\dbE[X(s)]\big\},\qq \dbE[u(s)]-\D(s)\dbE[X(s)]=v(s).$$
Then \rf{12.30} yields
$$\ba{lll}
0\4n&\les\4n&\ds\int_t^T\blan\bar\Si(s)\big\{\dbE[u(s)]-\D(s)\dbE[X(s)]\big\},\dbE[u(s)]-\D(s)\dbE[X(s)]\bran ds\\
\ns\4n&=\4n&\ds \int_t^T\blan\bar\Si(s)v(s),v(s)\bran ds,
\qq \forall v(\cd)\in L^2(t,T;\dbR^m),\ea$$
which implies that $\bar\Si\ges0$. The proof is completed.
\endpf

\section{Characterization of closed-loop solvability}

The aim of this section is to provide a characterization of the closed-loop solvability
of Problem (MF-LQ) in terms of the GREs \rf{Ric}, a linear backward stochastic differential
equation (BSDE, for short), and a linear terminal value problem of ODE.
In the case of Problem (MF-LQ)$^0$, it turns out that the regular solvability
of the GREs \rf{Ric} is not only necessary but also sufficient for the existence
of an optimal closed-loop strategy.

\bt{S&N}\sl Let {\rm(H1)--(H2)} hold and $t\in(0,T)$. Then Problem {\rm(MF-LQ)}
is closed-loop solvable on $[t,T]$ if and only if the following hold:
\begin{enumerate}[~~\,\rm(i)]
\item The GREs \rf{Ric} admits a regular solution
$(P(\cd),\Pi(\cd))\in C([t,T];\dbS^n)\times C([t,T];\dbS^n)$.
\item The adapted solution $(\eta(\cd),\z(\cd))$ to the BSDE
\bel{eta-zeta}\left\{\2n\ba{ll}
\ds d\eta(s)=-\big[(A+B\Th)^\top\eta+(C+D\Th)^\top\z+(C+D\Th)^\top P\si\\
\ns\ds\qq\qq~~~+\Th^\top\rho+Pb+q\big]ds+\z dW(s),\qq s\in[t,T],\\
\ns\ds\eta(T)=g,\ea\right.\ee
satisfies
$$\left\{\2n\ba{ll}
\ds B^\top(\eta-\dbE[\eta])+D^\top(\z-\dbE[\z])+D^\top P(\si-\dbE[\si])+\rho-\dbE[\rho]\in\cR(\Si),
\qq\ae~s\in[t,T],~\as\\
\ns\ds \f\equiv-\Si^\dag\big\{B^\top(\eta-\dbE[\eta])+D^\top(\z-\dbE[\z])+D^\top P(\si-\dbE[\si])
+\rho-\dbE[\rho]\big\}\in L_\dbF^2(t,T;\dbR^m),
\ea\right.$$
and the solution $\bar\eta(\cd)$ to the ODE
$$\left\{\2n\ba{ll}
\ds\dot{\bar\eta}+\big[(A+\bar A)+(B+\bar B)\G\big]^\top\bar\eta
+\G^\top\Big\{(D+\bar D)^\top\big(P\dbE[\si]+\dbE[\z]\big)+\dbE[\rho]+\bar\rho\Big\}\\
\ns\ds~~\1n+(C+\bar C)^\top\big(P\dbE[\si]+\dbE[\z]\big)+\dbE[q]+\bar q+\Pi\dbE[b]=0,\qq\ae~s\in[t,T],\\
\ns\ds\bar\eta(T)=\dbE[g]+\bar g,\ea\right.$$
satisfies
$$\left\{\2n\ba{ll}
\ds (B+\bar B)^\top\bar\eta+(D+\bar D)^\top(P\dbE[\si]+\dbE[\z])
+\dbE[\rho]+\bar\rho\in\cR(\bar\Si),\qq\ae~s\in[t,T],\\
\ns\ds\bar\f\equiv-\bar\Si^\dag\big\{(B+\bar B)^\top\bar\eta
+(D+\bar D)^\top\big(P\dbE[\si]+\dbE[\z]\big)+\dbE[\rho]+\bar\rho\big\}\in L^2(t,T;\dbR^m),
\ea\right.$$
where
$$\left\{\2n\ba{lll}
\ds \Si\4n&=\4n&\ds R+D^\top PD,\\
\ns\ds \bar\Si\4n&=\4n&\ds R+\bar R+(D+\bar D)^\top P(D+\bar D),\\
\ns\ds \Th\4n&=\4n&\ds-\,\Si^\dag(B^\top P+D^\top PC+S),\\
\ns\ds \G\4n&=\4n&\ds-\,\bar\Si^\dag\big[(B+\bar B)^\top\Pi+(D+\bar D)^\top P(C+\bar C)+(S+\bar S)\big].
\ea\right.$$
\end{enumerate}

In the above case, the optimal closed-loop strategy $(\Th^*(\cd),\bar\Th^*(\cd),u^*(\cd))$
admits the following representation:
\bel{close-rep}\left\{\2n\ba{lll}
\ds \Th^*\4n&=\4n&\ds \Th+(I-\Si^\dag\Si)\th,\\
\ns\ds\bar\Th^*\4n&=\4n&\ds \G-\Th+(I-\bar\Si^\dag\bar\Si)\t-(I-\Si^\dag\Si)\th,\\
\ns\ds u^*\4n&=\4n&\ds \f+\bar\f+(I-\Si^\dag\Si)(\nu-\dbE[\nu])+(I-\bar\Si^\dag\bar\Si)\bar\nu,
\ea\right.\ee
where $\th(\cd),\t(\cd)\in L^2(t,T;\dbR^{m\times n})$, $\nu(\cd)\in L_\dbF^2(t,T;\dbR^m)$,
and $\bar\nu(\cd)\in L^2(t,T;\dbR^m)$. Moreover, the value $V(t,\xi)$ is given by
\bel{V-rep}\ba{ll}
\ds V(t,\xi)=\dbE\lan P(t)(\xi-\dbE[\xi])+2\eta(t),\xi-\dbE[\xi]\ran
+\lan\Pi(t)\dbE[\xi]+2\bar\eta(t),\dbE[\xi]\ran\\
\ns\ds\qq\q\ ~~
+\dbE\int_t^T\Big\{\lan P\si,\si\ran+2\lan\eta,b-\dbE[b]\ran
+2\lan\z,\si\ran+2\lan\bar\eta,\dbE[b]\ran
-\lan\Si\f,\f\ran-\blan\bar\Si\bar\f,\bar\f\bran\Big\}ds.\ea\ee
\et

\it Proof. \rm {\it Necessity.} Suppose that $(\Th^*(\cd),\bar\Th^*(\cd),u^*(\cd))\in\sC[t,T]$
is an optimal closed-loop strategy of Problem (MF-LQ) on $[t,T]$. Then it follows from
Theorem \ref{Th-Nec} that the GREs \rf{Ric} admits a regular solution $(P(\cd),\Pi(\cd))$.
To determine $u^*(\cd)$, let $(X^*(\cd),Y^*(\cd),Z^*(\cd))$ be the adapted solution of
\rf{MF-FBSDE}. Proceeding as in the proof of Proposition \ref{prop2}, we see that
$(X^*(\cd),Y^*(\cd),Z^*(\cd))$ satisfies \rf{yueshu}.
Now, let $\D(\cd)=\Th^*(\cd)+\bar \Th^*(\cd)$ and define
\bel{a-b}\left\{\2n\ba{lll}
\ds\a\4n&=\4n&\ds(Y^*-\dbE[Y^*])-P(X^*-\dbE[X^*]),\\
\ns\ds\b\4n&=\4n&\ds Z^*-P(C+D\Th^*)(X^*-\dbE[X^*])-PD(u^*-\dbE[u^*])-P\si\\
\4n&~\4n&\ds-\,P\big[C+\bar C+(D+\bar D)\D\big]\dbE[X^*]-P(D+\bar D)\dbE[u^*],\\
\ns\ds\wt\eta\4n&=\4n&\ds\dbE[Y^*]-\Pi\dbE[X^*].
\ea\right.\ee
We have the following:
\bel{Eb}\dbE[\b]=\dbE[Z^*]-P\dbE[\si]-P\big[C+\bar C+(D+\bar D)\D\big]\dbE[X^*]-P(D+\bar D)\dbE[u^*],\ee\\[-3.5em]
\bel{b-Eb}\b-\dbE[\b]=Z^*-\dbE[Z^*]-P(C+D\Th^*)\big(X^*-\dbE[X^*]\big)-PD\big(u^*-\dbE[u^*]\big)-P(\si-\dbE[\si]).\ee
Recall from the proof of Theorem \ref{Th-Nec} that
\bel{10.1-P}\left\{\2n\ba{ll}
\ds \dot P+(A+B\Th^*)^\top P+P(A+B\Th^*)+(C+D\Th^*)^\top P(C+D\Th^*)\\
\ns\ds\q+\,Q+(\Th^*)^\top S+S^\top\Th^*+(\Th^*)^\top R\Th^*=0,\\
\ns\ds\Si\Th^*+B^\top P+D^\top PC+S=0,\\
\ns\ds \Th^*=\Th+(I-\Si^\dag\Si)\th \q\hb{for some }\th(\cd)\in L^2(t,T;\dbR^{m\times n}),
\ea\right.\ee
and
\bel{10.1-Pi}\left\{\2n\ba{ll}
\ds \dot\Pi+\big[A+\bar A+(B+\bar B)\D\big]^\top\Pi+\Pi\big[A+\bar A+(B+\bar B)\D\big]\\
\ns\ds\q+\,Q+\bar Q+(C+\bar C)^\top P(C+\bar C)+\D^\top\big[R+\bar R+(D+\bar D)^\top P(D+\bar D)\big]\D\\
\ns\ds\q+\,\D^\top\big[(D+\bar D)^\top P(C+\bar C)+S+\bar S\big]
+\big[(D+\bar D)^\top P(C+\bar C)+S+\bar S\big]^\top\D=0,\\
\ns\ds\bar\Si\D+(B+\bar B)^\top\Pi+(D+\bar D)^\top P(C+\bar C)+S+\bar S=0,\\
\ns\ds \D=\G+\big(I-\bar\Si^\dag\bar\Si\big)\t \q\hb{for some }\t(\cd)\in L^2(t,T;\dbR^{m\times n}).
\ea\right.\ee
Then we have $\a(T)=g-\dbE[g]$ and
$$\ba{lll}
\ds d\a\4n&=\4n&\ds d(Y^*-\dbE[Y^*])-\dot P(X^*-\dbE[X^*])ds-Pd(X^*-\dbE[X^*])\\
\ns\4n&=\4n&\ds-\,\Big\{(A+B\Th^*)^\top(Y^*-\dbE[Y^*])+(C+D\Th^*)^\top(Z^*-\dbE[Z^*])\\
\ns\4n&~\4n&\ds\q~~+\wt Q(X^*-\dbE[X^*])+\wt S^\top(u^*-\dbE[u^*])+\wt q-\dbE[\wt q\,]\Big\}ds+Z^*dW\\
\ns\4n&~\4n&\ds-\,\dot P(X^*-\dbE[X^*])ds-P\Big\{(A+B\Th^*)(X^*-\dbE[X^*])+B(u^*-\dbE[u^*])+b-\dbE[b]\Big\}ds\\
\ns\4n&~\4n&\ds-\,P\Big\{(C+D\Th^*)(X^*-\dbE[X^*])+\big[C+\bar C+(D+\bar D)\D\big]\dbE[X^*]\\
\ns\4n&~\4n&\ds\qq~+D(u^*-\dbE[u^*])+(D+\bar D)\dbE[u^*]+\si\Big\}dW\\
\ns\4n&=\4n&\ds-\,\Big\{(A+B\Th^*)^\top\a+(A+B\Th^*)^\top P(X^*-\dbE[X^*])\\
\ns\4n&~\4n&\ds\q~~+(C\1n+\1nD\Th^*)^\top\big[\b\1n-\1n\dbE[\b]\1n
+\1nP(C\1n+\1nD\Th^*)(X^*\1n-\1n\dbE[X^*])\1n+\1nPD(u^*\1n-\1n\dbE[u^*])\1n+\1nP(\si\1n-\1n\dbE[\si])\big]\\
\ns\4n&~\4n&\ds\q~~+\wt Q(X^*-\dbE[X^*])+\big[S^\top+(\Th^*)^\top R\big](u^*-\dbE[u^*])
+q-\dbE[q]+(\Th^*)^\top(\rho-\dbE[\rho])\\
\ns\4n&~\4n&\ds\q~~+\dot P(X^*\1n-\1n\dbE[X^*])\1n+\1nP(A\1n+\1nB\Th^*)(X^*\1n-\1n\dbE[X^*])\1n
+\1nPB(u^*\1n-\1n\dbE[u^*])\1n+\1nP(b\1n-\1n\dbE[b])\Big\}ds\1n+\1n\b dW\\
\ns\4n&=\4n&\ds-\,\Big\{(A+B\Th^*)^\top\a+(C+D\Th^*)^\top(\b-\dbE[\b])+(C+D\Th^*)^\top P(\si-\dbE[\si])\\
\ns\4n&~\4n&\ds\q~~+(\Th^*)^\top(\rho-\dbE[\rho])+P(b-\dbE[b])+q-\dbE[q]\Big\}ds+\b dW\\
\ns\4n&~\4n&\ds-\,\Big\{\big[\dot P+(A+B\Th^*)^\top P+P(A+B\Th^*)
+(C+D\Th^*)^\top P(C+D\Th^*)+\wt Q\,\big](X^*-\dbE[X^*])\\
\ns\4n&~\4n&\ds\q~~+\big[(C+D\Th^*)^\top PD+S^\top+(\Th^*)^\top R+PB\big](u^*-\dbE[u^*])\Big\}ds\\
\ns\4n&=\4n&\ds-\,\Big\{(A+B\Th^*)^\top\a+(C+D\Th^*)^\top(\b-\dbE[\b])+(C+D\Th^*)^\top P(\si-\dbE[\si])\\
\ns\4n&~\4n&\ds\q~~+(\Th^*)^\top(\rho-\dbE[\rho])+P(b-\dbE[b])+q-\dbE[q]\Big\}ds+\b dW.
\ea$$
Also, we have $\wt\eta(T)=\dbE[g]+\bar g$ and
$$\ba{lll}
\ds-{d\wt\eta\over ds}\4n&=\4n&\ds-\,{d\dbE[Y^*]\over ds}+\dot\Pi\dbE[X^*]+\Pi{d\dbE[X^*]\over ds}\\
\ns\4n&=\4n&\ds\big[A+\bar A+(B+\bar B)\D\big]^\top\dbE[Y^*]+\big[C+\bar C+(D+\bar D)\D\big]^\top\dbE[Z^*]\\
\ns\4n&~\4n&\ds+\,\big(\wt Q+\h Q\big)\dbE[X^*]+\big(\wt S+\h S\,\big)^\top\dbE[u^*]
+\dbE[\wt q\,]+\h q+\dot\Pi\dbE[X^*]\\
\ns\4n&~\4n&\ds+\,\Pi\big[A+\bar A+(B+\bar B)\D\big]\dbE[X^*]+\Pi(B+\bar B)\dbE[u^*]+\Pi\dbE[b]\\
\ns\4n&=\4n&\ds\big[A+\bar A+(B+\bar B)\D\big]^\top\wt\eta+\big[A+\bar A+(B+\bar B)\D\big]^\top\Pi\dbE[X^*]\\
\ns\4n&~\4n&\ds+\,\big[C\1n+\1n\bar C\1n+\1n(D\1n+\1n\bar D)\D\big]^\top\Big\{\dbE[\b]\1n+\1nP\dbE[\si]\1n
+\1nP\big[C\1n+\1n\bar C\1n+\1n(D\1n+\1n\bar D)\D\big]\dbE[X^*]\1n+\1nP(D\1n+\1n\bar D)\dbE[u^*]\Big\}\\
\ns\4n&~\4n&\ds+\,\big(\wt Q+\h Q\big)\dbE[X^*]+\big(\wt S+\h S\,\big)^\top\dbE[u^*]
+\dbE[q]+\bar q+\D^\top(\dbE[\rho]+\bar\rho)+\dot\Pi\dbE[X^*]\\
\ns\4n&~\4n&\ds+\,\Pi\big[A+\bar A+(B+\bar B)\D\big]\dbE[X^*]+\Pi(B+\bar B)\dbE[u^*]+\Pi\dbE[b]\\
\ns\4n&=\4n&\ds\big[A+\bar A+(B+\bar B)\D\big]^\top\wt\eta
+\D^\top\Big\{(D+\bar D)^\top(P\dbE[\si]+\dbE[\b])+\dbE[\rho]+\bar\rho\Big\}\\
\ns\4n&~\4n&\ds+\,(C+\bar C)^\top(P\dbE[\si]+\dbE[\b])+\dbE[q]+\bar q+\Pi\dbE[b]\\
\ns\4n&~\4n&\ds+\,\Big\{\dot\Pi+\big[A+\bar A+(B+\bar B)\D\big]^\top\Pi+\Pi\big[A+\bar A+(B+\bar B)\D\big]\\
\ns\4n&~\4n&\ds\q~~+\big[C+\bar C+(D+\bar D)\D\big]^\top P\big[C+\bar C+(D+\bar D)\D\big]+\wt Q+\h Q\Big\}\dbE[X^*]\\
\ns\4n&~\4n&\ds+\,\Big\{\Pi(B+\bar B)+\big[C+\bar C+(D+\bar D)\D\big]^\top P(D+\bar D)
+\big(\wt S+\h S\,\big)^\top\Big\}\dbE[u^*]\\
\ns\4n&=\4n&\ds\big[A+\bar A+(B+\bar B)\D\big]^\top\wt\eta
+\D^\top\Big\{(D+\bar D)^\top(P\dbE[\si]+\dbE[\b])+\dbE[\rho]+\bar\rho\Big\}\\
\ns\4n&~\4n&\ds+\,(C+\bar C)^\top(P\dbE[\si]+\dbE[\b])+\dbE[q]+\bar q+\Pi\dbE[b].
\ea$$
Moreover, we have from \rf{yueshu}:
\bel{3.7-1}(R+\bar R)\dbE[u^*]+(B+\bar B)^\top\dbE[Y^*]+(D+\bar D)^\top\dbE[Z^*]
+(\wt S+\h S\,)\dbE[X^*]+\dbE[\rho]+\bar\rho=0,\ee\\[-3.5em]
\bel{3.7-2}R(u^*-\dbE[u^*])+ B^\top(Y^*-\dbE[Y^*])+D^\top(Z^*-\dbE[Z^*])
+\wt S(X^*-\dbE[X^*])+\rho-\dbE[\rho]=0.\ee
Now \rf{a-b}, \rf{Eb}, \rf{10.1-Pi} and \rf{3.7-1} yield
$$\ba{lll}
\ds 0\4n&=\4n&\ds (R+\bar R)\dbE[u^*]+(B+\bar B)^\top\wt\eta+(B+\bar B)^\top\Pi\dbE[X^*]\\
\ns\4n&~\4n&\ds+\,(D+\bar D)^\top(\dbE[\b]+P\dbE[\si])+(D+\bar D)^\top P\big[C+\bar C+(D+\bar D)\D\big]\dbE[X^*]\\
\ns\4n&~\4n&\ds+\,(D+\bar D)^\top P(D+\bar D)\dbE[u^*]+(\wt S+\h S\,)\dbE[X^*]+\dbE[\rho]+\bar\rho\\
\ns\4n&=\4n&\ds\bar\Si\dbE[u^*]+(B+\bar B)^\top\wt\eta+(D+\bar D)^\top(\dbE[\b]+P\dbE[\si])+\dbE[\rho]+\bar\rho.
\ea$$
Hence,
$$(B+\bar B)^\top\wt\eta+(D+\bar D)^\top(\dbE[\b]+P\dbE[\si])+\dbE[\rho]+\bar\rho\in\cR(\bar\Si).$$
Since $\bar\Si^\dag\big[(B+\bar B)^\top\wt\eta+(D+\bar D)^\top(\dbE[\b]+P\dbE[\si])+\dbE[\rho]+\bar\rho\big]
=-\bar\Si^\dag\bar\Si\dbE[u^*]$ and $\bar\Si^\dag\bar\Si$ is an orthogonal projection, we have
$$\bar\Si^\dag\big\{(B+\bar B)^\top\wt\eta+(D+\bar D)^\top(\dbE[\b]+P\dbE[\si])
+\dbE[\rho]+\bar\rho\big\}\in L^2(t,T;\dbR^m),$$\\[-3.5em]
\bel{Eu-star}\dbE[u^*]=-\bar\Si^\dag\big\{(B+\bar B)^\top\wt\eta+(D+\bar D)^\top(\dbE[\b]+P\dbE[\si])
+\dbE[\rho]+\bar\rho\big\}+(I-\bar\Si^\dag\bar\Si)\bar\nu,\ee
for some $\bar\nu(\cd)\in L^2(t,T;\dbR^m)$. Consequently,
$$\ba{lll}
\ds-\dot{\wt\eta}\4n&=\4n&\ds(A+\bar A)^\top\wt\eta
+\D^\top\Big\{(B+\bar B)^\top\wt\eta+(D+\bar D)^\top(P\dbE[\si]+\dbE[\b])+\dbE[\rho]+\bar\rho\Big\}\\
\ns\4n&~\4n&\ds+\,(C+\bar C)^\top(P\dbE[\si]+\dbE[\b])+\dbE[q]+\bar q+\Pi\dbE[b]\\
\ns\4n&=\4n&\ds(A+\bar A)^\top\wt\eta+\big[\G+(I-\bar\Si^\dag\bar\Si)\t\big]^\top
\Big\{(B+\bar B)^\top\wt\eta+(D+\bar D)^\top(P\dbE[\si]+\dbE[\b])+\dbE[\rho]+\bar\rho\Big\}\\
\ns\4n&~\4n&\ds+\,(C+\bar C)^\top(P\dbE[\si]+\dbE[\b])+\dbE[q]+\bar q+\Pi\dbE[b]\\
\ns\4n&=\4n&\ds\big[A+\bar A+(B+\bar B)\G\big]^\top\wt\eta
+\G^\top\Big\{(D+\bar D)^\top(P\dbE[\si]+\dbE[\b])+\dbE[\rho]+\bar\rho\Big\}\\
\ns\4n&~\4n&\ds+\,(C+\bar C)^\top(P\dbE[\si]+\dbE[\b])+\dbE[q]+\bar q+\Pi\dbE[b].
\ea$$
Likewise, \rf{a-b}, \rf{b-Eb}, \rf{10.1-P} and \rf{3.7-2} yield
$$\Si(u^*-\dbE[u^*])+B^\top\a+D^\top(\b-\dbE[\b])+D^\top P(\si-\dbE[\si])+\rho-\dbE[\rho]=0,$$
which implies that
$$B^\top\a+D^\top(\b-\dbE[\b])+D^\top P(\si-\dbE[\si])+\rho-\dbE[\rho]\in\cR(\Si),$$
$$\Si^\dag\big\{B^\top\a+D^\top(\b-\dbE[\b])+D^\top P(\si-\dbE[\si])+\rho-\dbE[\rho]\big\}
\in L_\dbF^2(t,T;\dbR^m),$$
and that
\bel{u-Eu-star}u^*-\dbE[u^*]=
-\Si^\dag\big\{B^\top\a+D^\top(\b-\dbE[\b])+D^\top P(\si-\dbE[\si])+\rho-\dbE[\rho]\big\}
+(I-\Si^\dag\Si)(\nu-\dbE[\nu]),\ee
for some $\nu(\cd)\in L_\dbF^2(t,T;\dbR^m)$. Consequently,
$$\ba{lll}
\ds d\a\4n&=\4n&\ds-\,\Big\{A^\top\a+C^\top(\b-\dbE[\b])+C^\top P(\si-\dbE[\si])+P(b-\dbE[b])+q-\dbE[q]\\
\ns\4n&~\4n&\ds\q~~+(\Th^*)^\top\big(B^\top\a+D^\top(\b-\dbE[\b])+D^\top P(\si-\dbE[\si])+\rho-\dbE[\rho]\big)\Big\}ds+\b dW\\
\ns\4n&=\4n&\ds-\,\Big\{A^\top\a+C^\top(\b-\dbE[\b])+C^\top P(\si-\dbE[\si])+P(b-\dbE[b])+q-\dbE[q]\\
\ns\4n&~\4n&\ds\q~~+\big[\Th+(I-\Si^\dag\Si)\th\big]^\top\big(B^\top\a+D^\top(\b-\dbE[\b])+D^\top P(\si-\dbE[\si])+\rho-\dbE[\rho]\big)\Big\}ds+\b dW\\
\ns\4n&=\4n&\ds-\,\Big\{(A+B\Th)^\top\a+(C+D\Th)^\top(\b-\dbE[\b])+(C+D\Th)^\top P(\si-\dbE[\si])\\
\ns\4n&~\4n&\ds\q~~+\Th^\top(\rho-\dbE[\rho])+P(b-\dbE[b])+q-\dbE[q]\Big\}ds+\b dW.
\ea$$
Since the adapted solution $(\eta(\cd),\z(\cd))$ of \rf{eta-zeta} satisfies
$$\left\{\2n\ba{ll}
\ds d(\eta-\dbE[\eta])=-\Big\{(A+B\Th)^\top(\eta-\dbE[\eta])+(C+D\Th)^\top(\z-\dbE[\z])
+(C+D\Th)^\top P(\si-\dbE[\si])\\
\ns\ds\qq\qq\qq\q~+\Th^\top(\rho-\dbE[\rho])+P(b-\dbE[b])+q-\dbE[q]\Big\}ds+\z dW(s),\qq s\in[t,T],\\
\ns\ds\eta(T)-\dbE[\eta(T)]=g-\dbE[g],\ea\right.$$
by the uniqueness of solutions, we must have $\a(\cd)=\eta(\cd)-\dbE[\eta(\cd)]$ and $\b(\cd)=\z(\cd)$,
and hence also $\wt\eta(\cd)=\bar\eta(\cd)$. Then, \rf{Eu-star} and \rf{u-Eu-star} yield
$$u^*=\f+(I-\Si^\dag\Si)(\nu-\dbE[\nu])+\bar\f+(I-\bar\Si^\dag\bar\Si)\bar\nu.$$
This proves the necessity, as well as \rf{close-rep}.

\ms

\it Sufficiency. \rm The proof is much like that of \cite[Theorem 5.2]{Sun}.
Let $(\Th^*(\cd),\bar\Th^*(\cd),u^*(\cd))$ be defined by \rf{close-rep}.
Then we have
\bel{10.3-1}\left\{\2n\ba{ll}
\ds B^\top P+D^\top PC+S=-\Si\Th=-\Si\Th^*,\\
\ns\ds B^\top(\eta-\dbE[\eta])+D^\top(\z-\dbE[\z])+D^\top P(\si-\dbE[\si])+\rho-\dbE[\rho]=-\Si\f,\\
\ns\ds(B+\bar B)^\top\Pi+(D+\bar D)^\top P(C+\bar C)+(S+\bar S)=-\bar\Si\G=-\bar\Si(\Th^*+\bar\Th^*),\\
\ns\ds(B+\bar B)^\top\bar\eta+(D+\bar D)^\top(P\dbE[\si]+\dbE[\z])+\dbE[\rho]+\bar\rho=-\bar\Si\bar\f.
\ea\right.\ee
For any $\xi\in L^2_{\cF_t}(\Om;\dbR^n)$ and $u(\cd)\in\cU[t,T]$, let $X(\cd)\equiv X(\cd\,;t,\xi,u(\cd))$
be the corresponding solution of \rf{state}. Proceeding similarly to the proof of \cite[Theorem 5.2]{Sun}
and using \rf{10.3-1}, we obtain
\bel{10.4-1}\ba{ll}
\ds J(t,\xi;u(\cd))-\dbE\lan P(t)(\xi-\dbE[\xi])+2\eta(t),\xi-\dbE[\xi]\ran
-\lan\Pi(t)\dbE[\xi]+2\bar\eta(t),\dbE[\xi]\ran\\
\ns\ds=\dbE\int_t^T\Big\{\lan P\si,\si\ran+2\lan\eta,b-\dbE[b]\ran
+2\lan\z,\si\ran+2\lan\bar\eta,\dbE[b]\ran-\lan\Si\f,\f\ran-\blan\bar\Si\bar\f,\bar\f\bran\\
\ns\ds\qq\qq~+\blan\Si\big\{u-\dbE[u]-\Th(X-\dbE[X])-\f\big\},u-\dbE[u]-\Th(X-\dbE[X])-\f\bran\\
\ns\ds\qq\qq~+\blan\bar\Si\big(\dbE[u]-\G\dbE[X]-\bar\f\big),\dbE[u]-\G\dbE[X]-\bar\f\bran\Big\}ds.
\ea\ee
Since $\Si,\bar\Si\ges0$, \rf{10.4-1} implies that
$$\ba{lll}
\ds J(t,\xi;u(\cd))\4n&\ges\4n&\ds\dbE\lan P(t)(\xi-\dbE[\xi])+2\eta(t),\xi-\dbE[\xi]\ran
+\lan\Pi(t)\dbE[\xi]+2\bar\eta(t),\dbE[\xi]\ran\\
\ns\4n&~\4n&\ds+\,\dbE\int_t^T\Big\{\lan P\si,\si\ran+2\lan\eta,b-\dbE[b]\ran
+2\lan\z,\si\ran+2\lan\bar\eta,\dbE[b]\ran-\lan\Si\f,\f\ran-\blan\bar\Si\bar\f,\bar\f\bran\Big\}ds\\
\ns\4n&=\4n&\ds J(t,\xi;\Th^*(\cd)X^*(\cd)+\bar\Th^*(\cd)\dbE[X^*(\cd)]+u^*(\cd)),
\q\1n~\forall(\xi,u(\cd))\in L^2_{\cF_t}(\Om;\dbR^n)\times\cU[t,T].
\ea$$
Therefore, $(\Th^*(\cd),\bar\Th^*(\cd),u^*(\cd))$ is an optimal closed-loop strategy
of Problem (MF-LQ) on $[t,T]$ and \rf{V-rep} holds. The proof is completed.
\endpf

\ms

In the special case that $b(\cd)$, $\si(\cd)$, $g$, $\bar g$, $q(\cd)$, $\bar q(\cd)$, $\rho(\cd)$,
and $\bar\rho(\cd)$ vanish, if the GREs \rf{Ric} admits a regular solution $(P(\cd),\Pi(\cd))$,
then condition (ii) of Theorem \ref{S&N} holds automatically. Indeed, one can easily check that
$(\eta(\cd),\z(\cd))=(0,0)$ and $\bar\eta(\cd)=0$. Thus, we have the following corollary.

\bc{}\sl Let {\rm(H1)--(H2)} hold and $t\in(0,T)$. Then Problem {\rm(MF-LQ)$^0$} is
closed-loop solvable on $[t,T]$ if and only if the GREs \rf{Ric} is regularly solvable.
\ec

\section{Conclusion}

This is an important yet challenging research topic. Recently there has been increasing interest
in studying this type of stochastic control problems as well as their applications. Beside this
work, the optimal stochastic control problems under MF-SDEs are underdeveloped in the literature,
and therefore many fundamental questions remain open and methodologies need to be significantly
improved. To establish new theory and hopefully to shed light on financial investment, we expect
the findings of this research program to add to various streams of the literature, such as portfolio
selection, optimal control techniques, financial risk management, and relative performance evaluation.


\begin{thebibliography}{90}
\addtolength{\itemsep}{-1.0ex}


\bibitem{Ait Rami-Moore-Zhou 2001} M.~Ait~Rami, J.~B.~Moore, and X.~Y.~Zhou,
    \it Indefinite stochastic linear quadratic control and generalized differential Riccati equation,
    \sl SIAM J. Control Optim., \rm 40 (2001), 1296--1311.

\bibitem{Andersson-Djehiche 2011} D.~Andersson and B.~Djehiche,
    \it A Maximum Principle for SDEs of Mean-Field Type,
    \sl Appl. Math. Optim., \rm 63 (2011), 341--356.

\bibitem{Athans 1968} M.~Athans,
    \it The matrix minimum principle,
    \sl Inform. and Control, \rm 11 (1968), 592--606.

\bibitem{Buckdahn-Djehiche-Li 2011} R.~Buckdahn, B.~Djehiche, and J.~Li,
    \it A general stochastic maximum principle for SDEs of mean-field type,
    \sl Appl. Math. Optim., \rm 64 (2011), 197--216.

\bibitem{Buckdahn-Djehiche-Li-Peng 2009} R.~Buckdahn, B.~Djehiche, J.~Li, and S.~Peng,
    \it Mean-field backward stochastic differential equations: a limit approach,
    \sl Ann. Probab., \rm 37 (2009), 1524--1565.

\bibitem{Buckdahn-Li-Peng 2009} R.~Buckdahn, J.~Li, and S.~Peng,
    \it Mean-field backward stochastic differential equations and related partial
        differential equations,
    \sl Stoch. Proc. Appl., \rm 119 (2009), 3133--3154.

\bibitem{Chen-Li-Zhou 1998} S.~Chen, X.~Li, and X.~Y.~Zhou,
    \it Stochastic linear quadratic regulators with indefinite control weight costs,
    \sl SIAM J. Control Optim., \rm 36 (1998), 1685--1702.

\bibitem{Chen-Yong 2000} S.~Chen and J.~Yong,
    \it Stochastic linear quadratic optimal control problems with random coefficients,
    \sl Chin. Ann. Math., \rm 21 B (2000), 323--338.

\bibitem{Cui-Li-Li 2014} X.~Y.~Cui, X.~Li, and D.~Li,
    \it Unified framework of mean-field formulations for optimal multi-period
        mean-variance portfolio selection,
    \sl IEEE Trans. Auto. Control, \rm 59 (2014), 1833--1844.

\bibitem{Elliott-Li-Ni 2013} R.~Elliott, X.~Li, and Y.~H.~Ni,
    \it Discrete time mean-field stochastic linear-quadratic optimal control problems,
    \sl Automatica, \rm 49 (2013), 3222--3233.

\bibitem{Huang-Li-Wang 2015} J.~Huang, X.~Li, and T.~X.~Wang,
    \it Mean-field linear-quadratic-Gaussian (LQG) games for stochastic integral systems,
    \sl IEEE Trans. Auto. Control, \rm Accepted for publication, 2015.

\bibitem{Huang-Li-Yong 2015} J.~Huang, X.~Li, and J.~Yong,
    \it A linear-quadratic optimal control problem for mean-field stochastic differential
        equations in infinite horizon,
    \sl Mathematical Control \& Related Fields, \rm 5 (2015), 97--139.

\bibitem{Kac 1956} M.~Kac,
    \it Foundations of kinetic theory,
    \sl Proc. the Third Berkeley Symposium on Mathematical Statistics and Probability,
    \rm 3 (1956), 171--197.

\bibitem{Li-Zhou 2006} X.~Li and X.~Y.~Zhou,
    \it Continuous-time mean-variance efficiency: The 80\% rule,
    \sl Ann. Appl. Probab., \rm 16 (2006), 1751--1763.

\bibitem{Li-Zhou-Lim 2001} X.~Li, X.~Y.~Zhou, and A.~E.~B.~Lim,
    \it Dynamic mean-variance portfolio selection with no-shorting constraints,
    \sl SIAM J. Control Optim., \rm 40 (2001), 1540--1555.

\bibitem{McKean 1966} H.~P.~McKean,
    \it A class of Markov processes associated with nonlinear parabolic equations,
    \sl Proc. of the National Academy of Sciences of the United States of America,
    \rm 56 (1966), 1907--1911.

\bibitem{Meyer-Oksendal-Zhou 2011} T.~Meyer-Brandis, B.~{\O}ksendal, and X.~Y.~Zhou,
    \it A mean-field stochastic maximum principle via Malliavin calculus.
        A special issue for Mark Davis' Festschrift,
    to appear in \sl Stochastics, \rm 84 (2012), 643--666.

\bibitem{Penrose 1955} R.~Penrose,
    \it A generalized inverse of matrices,
    \sl Proc. Cambridge Philos. Soc., \rm 52 (1955), 17--19.

\bibitem{Sun} J.~Sun,
    \it Mean-Field Stochastic Linear Quadratic Optimal Control Problems:
        Open-Loop Solvabilities,
    \sl arXiv: 1509.02100v2. \rm

\bibitem{Sun-Li-Yong} J.~Sun, X.~Li, and J.~Yong,
    \it Open-Loop and Closed-Loop Solvabilities for Stochastic Linear Quadratic
        Optimal Control Problems,
    \sl arXiv: 1508.02163. \rm

\bibitem{Sun-Yong 2014} J.~Sun and J.~Yong,
    \it Linear Quadratic Stochastic Differential Games:
        Open-Loop and Closed-Loop Saddle Points,
    \sl SIAM J. Control Optim., \rm 52 (2014), 4082--4121.

\bibitem{Sun-Yong-Zhang 2016} J.~Sun, J.~Yong, and S.~Zhang, \it Linear quadratic stochastic two-person zero-sum differential games in an infinite horizon, \sl ESAIM COCV, \rm to appear.

\bibitem{Tang 2003} S.~Tang,
    \it General linear quadratic optimal stochastic control problems with random coefficients:
        linear stochastic Hamilton systems and backward stochastic Riccati equations,
    \sl SIAM J. Control Optim., \rm 42 (2003), 53--75.

\bibitem{Wonham 1968} W.~M.~Wonham,
    \it On a matrix Riccati equation of stochastic control,
    \sl SIAM J. Control Optim., \rm 6 (1968), 681--697.

\bibitem{Yong 2013} J.~Yong,
    \it Linear-Quadratic Optimal Control Problems for Mean-Field Stochastic
        Differential Equations,
    \sl SIAM J. Control Optim., \rm 51 (2013), 2809--2838.

\bibitem{Yong-Zhou 1999} J.~Yong and X.~Y.~Zhou,
    \sl Stochastic Controls: Hamiltonian Systems and HJB Equations,
    \rm Springer-Verlag, New York, 1999.

\bibitem{Zhou-Li 2000} X.~Y.~Zhou and D.~Li,
    \it Continuous-time mean-variance portfolio selection: A stochastic LQ framework,
    \sl Appl. Math. Optim., \rm 42 (2000), 19--33.

\end{thebibliography}
\end{document}